\documentclass[12pt]{amsart}
\usepackage{latexsym,amssymb,amsmath,amsthm,graphicx,color}
\usepackage[section]{placeins} 
\definecolor{red1}{rgb}{1,0.9,0.9} \definecolor{blue1}{rgb}{0.9,0.9,1} \definecolor{green1}{rgb}{0.9,1,0.9}
\definecolor{yellow1}{rgb}{1,1,0.9} \definecolor{yellow2}{rgb}{1,1,0.8}
    
\def\Dcal{\mathcal{D}}  \def\Ecal{\mathcal{E}} \def\Ccal{\mathcal{C}} \def\Gcal{\mathcal{G}} 
\def\Scal{\mathcal{S}} \def\Bcal{\mathcal{B}}    
\def\Xcal{\mathcal{X}}  
     
\def\question#1{ \vspace{2mm} \begin{center} \fcolorbox{green1}{green1}{ \parbox{11.2cm}{{\bf Question:} #1}} \vspace{2mm} \end{center} }
\def\conjecture#1{ \vspace{2mm} \begin{center} \fcolorbox{green1}{green1}{ \parbox{11.2cm}{{\bf Conjecture:} #1}} \vspace{2mm} \end{center} }
\def\resultlemma#1{ \vspace{2mm} \begin{center} \fcolorbox{yellow1}{yellow1}{ \parbox{11.2cm}{{\bf Lemma:} #1}} \vspace{2mm} \end{center} }
\def\resulttheorem#1{ \vspace{2mm} \begin{center} \fcolorbox{yellow1}{yellow1}{ \parbox{11.2cm}{{\bf Theorem:} #1}} \vspace{2mm} \end{center} }
\def\resultcorollary#1{ \vspace{2mm} \begin{center} \fcolorbox{yellow1}{yellow1}{ \parbox{11.2cm}{{\bf Corollary:} #1}} \vspace{2mm} \end{center} }
\def\definition#1{ \vspace{2mm} \begin{center} \fcolorbox{red1}{red1}{ \parbox{11.2cm}{{\bf Definition:} #1}} \vspace{2mm} \end{center} }
\title{Graphs with Eulerian unit spheres}
\author{Oliver Knill}
\date{January 11, 2015}
\address{ Department of Mathematics \\ Harvard University \\ Cambridge, MA, 02138 }

\begin{document}
\maketitle

\begin{abstract}
$d$-spheres are defined graph theoretically and inductively as the empty graph in dimension $d=-1$ and 
$d$-dimensional graphs for which all unit spheres $S(x)$ are $(d-1)$-spheres and such that for $d \geq 0$ 
the removal of one vertex renders the graph contractible. Eulerian $d$-spheres are geometric $d$-spheres 
which can be colored with $d+1$ colors. They are Eulerian graphs in the classical sense and for $d \geq 2$,
all unit spheres of an Eulerian sphere are Eulerian spheres. We prove here that $G$ is an Eulerian 
sphere if and only if the degrees of all the $(d-2)$-dimensional sub-simplices in $G$ are even. 
This generalizes a result of Kempe-Heawood for $d=2$ and is work related to the conjecture 
that all $d$-spheres have chromatic number $d+1$ or $d+2$ which is based on the geometric conjecture 
that every $d$-sphere can be embedded in an Eulerian $(d+1)$-sphere.
For $d=2$, such an embedding into an Eulerian 3-sphere would lead to a 
geometric proof of the 4 color theorem, allowing to see ``why 4 colors suffice". 
To achieve the goal of coloring a $d$-sphere $G$ with $d+2$ colors, we hope 
to embed it into a $(d+1)$-sphere and refine or thin out the later using special homotopy deformations
without touching the embedded sphere. Once rendered Eulerian and so $(d+2)$-colorable, it colors the
embedded graph $G$. In order to define the degree of a simplex, we introduce a notion of 
dual graph $\hat{H}$ of a subgraph $H$ in a general finite simple graph $G$. This leads to a natural sphere
bundle over the simplex graph. We look at geometric graphs which admit a unique geodesic flow: their unit spheres must be 
Eulerian. We define Platonic spheres graph theoretically as $d$-spheres for which all unit spheres $S(x)$ 
are graph isomorphic Platonic $(d-1)$-spheres.
Gauss-Bonnet allows a straightforward classification within graph theory independent of the classical 
Schl\"afli-Schoute-Coxeter classification: all spheres are Platonic for $d \leq 1$, the octahedron and 
icosahedron are the Platonic $2$-spheres, the sixteen and six-hundred cells are the Platonic $3$-spheres. 
The cross polytope is the unique Platonic $d$-sphere for $d>3$. It is Eulerian. 
\end{abstract}

\section{The geometry of graphs}

When studying geometric graphs in a purely graph theoretical manner, there is no need to
carrying around Euclidean spaces. This is analogue to look at groups detached from
concrete representations given by transformations in Euclidean space. This detachment from 
Euclidean structures can be pedagogically helpful: abstract group theory can be grasped early 
by students as symmetry groups or Rubik type puzzles which are intuitive. 
Kids without abstract training in algebra can solve these finitely
presented groups and do so often more successfully than adults with abstract mathematical training. 
Graphs are intuitive too, as they can be drawn. 
Place five dots on paper and connect them all: this is the hyper-tetrahedron, an example of a four
dimensional graph. Want to see a 3-dimensional sphere? Draw a regular 8-gon, then connect all 
diagonals except the main diagonals. You might not see the sphere structure yet, but it can be done
by looking at a unit sphere $S(x)$ of a vertex $x$. It is the subgraph formed by all points 
connected to $x$. If we draw this 
sphere, it can be represented as the regular hexagon, where all diagonals except the main diagonals 
are placed. Why is this a 2-sphere? Because after removing one point, the remaining graph can be collapsed 
and because the unit spheres $S(x)$ are spheres: each unit sphere is a $4$-gon, which is
an example of a $1$-sphere. How so? Because removing one point in this quadrangle produces an 
interval graph with $3$ vertices which can be collapsed and because every unit sphere in this graph is a 
$0$-sphere consisting of two disconnected points. Why is the graph consisting of two disconnected points a $0$-sphere: 
because removing one point has made the graph collapsed to one point
and because every unit sphere is empty. We just have to impose the initial induction assumption that the empty 
graph is the $(-1)$-sphere and we have defined inductively a concept of ``sphere" in arbitrary dimension $d$.
This construction illustrates how important dimension is: the dimension of a graph is inductively defined
as $1$ plus the average dimension of its unit spheres. An other important ingredient has entered in the definition of spheres:
the notion of homototopy is based on the notion of reduction steps and also inductive: a graph can be collapsed 
to a point, if there exists a vertex $x$ such that both its unit sphere as well as the graph obtained by 
removing the point can be collapsed to a point.
Our 4 dimensional simplex drawn before can be collapsed to a point because every unit sphere
is a 3-dimensional simplex and also because removing the vertex deforms the graph to a 3-dimensional simplex. 
So, by induction a complete $d$-simplex is collapsible if the $(d-1)$-simplex is collapsible. If as an induction
base, the one point graph is assumed to be collapsible, we are done. Why is a punctured
sphere collapsible? Take the octahedron; removing a single vertex produces the wheel graph with a 4-gon
boundary. The later is collapsible because we take away any outer vertex  we get the kite graph, two triangles
glued along an edge. And the kite graph is collapsible, because we can get from it a triangle after removing a vertex. 
The triangle is collapsible because it is a simplex. 
As in the continuum, one would have to distinguish collapsibility with contractibility. The later means homotopic
to a point, where homotopy is a chain of collapse or inverse collapsing steps and is not equivalent. There
are graphs like the dunce hat which first need to be enlarged before one can contract them to a point. Also this
is analogue to the continuum \cite{josellisknill}. But this subtlety is not relevant for spheres. \\

Why has this simple setup for geometric graphs not appeared earlier? Maybe because it would not have been right:
only since the Poincar\'e conjecture was proven, we know that spheres are homotopy spheres. Since it is a theorem
that topological $d$-spheres can be characterized as geometric objects of dimension $d$ which have the homotopy 
type of spheres, this is a valid way how to define spheres in graph theory.
An intuitive definition of dimension for graphs is to add $1$ to the average dimension of the
unit spheres. This definition is motivated by a classical Menger-Urysohn dimension but the later
assigns to graphs dimension $0$, because the topology generated by the usual distance metric
renders graphs $0$-dimensional: the topology on a graph is the discrete topology. One has to refer to pointless topology
to get a useful notion of homeomorphism \cite{KnillTopology}. 
Graphs are usually perceived as low dimensional objects: they are usually treated as one-dimensional 
simplicial complexes with zero and one-dimensional cohomology only, despite the fact that the higher 
dimensional nature is acknowledged in computer graphics or computational topology. We read for example in
\cite{coxeter}: ``a one-dimensional complex is a graph".
Without having to impose more structure, graphs already have a natural simplicial complex associated with them. 
This structure of complete subgraphs inside the graph 
allows also conveniently to compute cohomology: less than two dozen lines of computer 
algebra code are needed to compute a basis for any cohomology group of a general finite simple graph. Without
additional libraries. The key is to build a matrix $D$ and to find the kernels of the block matrices of 
$D^2$. For code see \cite{DiracKnill}.  \\

As early steps in geometry show, the notions of spheres and lines are pivotal for our understanding of a geometric
space. Both are based on distance: spheres are points with equal distance to a given point and lines are
curves which locally minimize distance. They can already be complicated in simple situations like a non-flat
surface, where spheres are wave fronts and lines are geodesics. Also spheres of codimension larger than 1 like a circle in space are 
important. We will see that in geometric graphs, graphs for which all unit spheres are graph theoretically defined abstract spheres, 
there is a natural way to get spheres as dual graphs to complete subgraphs. The notion of ``line" in geometry 
is superseded by the notion of ``geodesic", shortest connections between points. For general geometric
graphs, there is no natural geodesic flow. Already for spheres like the icosahedron, we have to tell
how to propagate a light ray through a disc. This requires to tell which ``in" direction goes to which ``out" direction.
There is a problem in general. Assume we hit cross road, where 5 streets meet. Which of the
two ``opposite streets" do we chose? We could throw a dice each time when hitting such an odd degree vertex
but that would render the geodesic flow a random walk. For two dimensional graphs there is a natural geodesic
flow if the graph has even degree at every point. In other words, $2$-dimensional graphs for which a geodesic flow
exists, have to be Eulerian. 
This already leads us to the topic studied here and it merges with an other topic in graph theory, 
the theory of colorings. The connection is that Eulerian spheres are exactly the spheres which can be colored
with 3 colors. But is pretty clear that in order to define a geodesic flow in higher dimensions which has the property
that is unique, we need the unit spheres to have a bit of symmetry. We need a fixed point free involution $T$ 
on the unit sphere $S(x)$, telling how to get from the incoming ray to the outgoing ray. 
This existence of a fixed point free involution 
is a weak form of ``projective" as if the graph is large enough, the quotient $G/T$ is then again a graph, 
a projective plane. In any case, weakly projective $d$-spheres are Eulerian in the sense that they can be colored 
with $d+1$ colors.  \\

When investigating the class of $d$-spheres which allow a minimal coloring with $d+1$ colors, it is 
important also to understand spheres of positive co-dimension. Examples are $(d-1)$-dimensional unit spheres or 
$(d-2)$-dimensional intersections of neighboring unit spheres. This naturally leads to a notion of duality, as 
the dual graph of a vertex is $S(x)$, and the dual graph of an edge $e=(a,b)$ is the intersection $S(a) \cap S(b)$, 
the dual of a triangle the intersection of three spheres etc.
In the case $d=3$, the dual sphere of an edge is a one-dimensional circular graph and its length is what we 
called edge degree. In the case $d=4$, the dual sphere of a triangle is a circular graph and its length is called
the degree of the triangle. In this case, degree is the number of hyper-tetrahedra ``hinging" at the triangle. 
When subdividing an edge, the degrees of the maximal simplices in the dual of the edge changes
parity. We are interested in these numbers because if in a $d$-sphere, all 
the $(d-2)$-dimensional simplices in $G$ are even, then the graph can be minimally colored. While we can not 
color it initially with $d+1$ colors, we aim to get there eventually using refinement or collapsing processes
\cite{knillgraphcoloring}. \\

To understand this in higher dimensions, it is necessary to look at spheres in graphs a bit more generally.
It turns out that we can define for any subgraph $H$ of a finite simple graph $G$ a $G$-dual graph $\hat{H}$. 
Unlike the dual graph $\hat{G}$ for the geometric graph $G$ itself which is defined by the intersection 
properties of its highest dimensional simplices, the $G$-dual graph depends on the host graph
$G$. If $H$ is a subgraph of $K$ and $K$ is a subgraph of $G$ then the $G$-dual graph $\hat{K}$ of $K$ 
is a subgraph of the $G$-dual graph $\hat{H}$ of $H$. As often in duality, the involution property does not kick in directly 
but it applies after applying the dual operation once. The reason is that for many subgraphs $H$, the dual graph will 
be empty so that the dual operation lock many subgraphs into the $\emptyset \leftrightarrow G$ cycle. This happens for
every subgraph $H$ of $G$ which has diameter larger than $2$ in $G$. What is important for us is that for
complete subgraphs $K$ of dimension $k$, the dual graph $\hat{K}$ is a sphere of dimension $d-k-1$.  \\

The notion of complementary duality is motivated by school geometry, where we construct a line 
perpendicular to a given line with the compass by intersecting spheres located at the end
points of an interval and end points of an intervals form a $0$- sphere, the two new intersection
points gave us a dual object. Doing the construction again at the dual line brings back the original $0$-sphere. 
The same construction can be done in graph theory. The ``compass" is now the ability to 
draw spheres. It is essentially the only tool at our disposition in the discrete. But sphere geometry 
is such a mighty instrument that we do not lose much. Actually, it is surprising how much of geometry 
goes over to graph theory, not only to geometric graphs, but to the full class of finite simple graphs. \\

There is a inductively defined class $\Scal_d$ of $d$-spheres, graphs which have all the properties of $d$-dimensional
spheres in the continuum. Geometric graphs in dimension $d$ are defined as graphs for which every unit sphere is a 
$(d-1)$-sphere.  Induction is important because also lower dimensional spheres are relevant in geometry. 
In graph theory, they can be constructed by intersecting unit spheres. For example,
the intersection of two neighboring unit spheres $S(a),S(b)$ in a $d$-dimensional graph produces a $(d-2)$-sphere because
$S(a) \cap S(b)$ is a unit sphere in the unit sphere $S(b)$ and because the definition of spheres has been 
done recursively. As we will see, if we take a $k$-dimensional simplex $K$ in a $d$-dimensional geometric
graph and intersect all these unit spheres centered at vertices in $K=K_{k+1}$, we obtain a sphere of 
dimension $d-k-1$.  \\

The idea to use higher dimensional geometry to color planar graphs has emerged in the late 70'ies in particular
through the work of Steve Fisk \cite{Fisk1977a}. 
The excitement and repercussions about the computer assisted proof of the 4-color theorem which was achieved
at just about the same time when \cite{Fisk1977a} was written, might have dampened a brisker development using geometric 
ideas even so it has not stopped. One difficulty with the Fisk approach as well as others following his work 
is that it based on more difficult definitions from topology, in particular the notion of 
``triangulation" which is complicated and which can carry surprises, especially in higher dimensions. 
Coloring questions of triangulations in 
higher dimensions in particular are expected to depend on the class of triangulations. We bypass the continuum 
and define what a ``geometric graph", a "sphere" or a "ball" is entirely within graph theory in a recursive 
way. Avoiding definitions from the continuum does not mean to ignore Riemannian geometry; almost all intuition from the 
continuum can be carried over almost verbatim from differential topology to graph theory if one lets go the usual 
assumption treating graphs as one-dimensional objects. And almost all definitions are motivated from continuum notions. 
And there are not only plenty of questions but plenty of entire research areas which have not yet been carried over
from the continuum to this combinatorial setting \cite{KnillBaltimore}. \\

The fresh start with clear and simple definitions allows for bolder conjectures in \cite{knillgraphcoloring} like

\conjecture{
All $G \in \Gcal_2$ can be colored with $3$,$4$ or $5$ colors.
}

\conjecture{
Every $G \in \Scal_d$ is $d+1$ or $(d+2)$-colorable.
}

The later is settled for $d=2$, where it is known to be equivalent to the 4-color theorem \cite{knillgraphcoloring}.
The reason why these statements are likely to hold true is because of the conjectured
picture that every orientable $2$-dimensional geometric graph can be embedded in an Eulerian three
sphere bounding a four dimensional ball of dimension $4$, which has a simply connected interior
(in the non-orientable case, a Moebius turn needs to to leave into the 4-ball but this does not affect its
simply connectivity) and that every $d$-sphere can be embedded in an Eulerian $(d+1)$-sphere.  
In order to investigate this further, we have to understand Eulerian spheres in 
any dimension and understand a special class of homotopy transformations called edge
subdivision or edge collapses. We would eventually like to know whether any two spheres 
can be deformed into each other by such homotopies and whether we can perform this talk under the 
``handicap" of not modifing an embedded smaller dimensional surface. 

\section{Graphs with Eulerian spheres}

The definition of geometric graphs and spheres is by induction, starting with the assumption 
that the empty graph $\emptyset$ is the $(-1)$-dimensional sphere so that $\Gcal_{-1}=\Scal_{-1}=\{\emptyset \}$
and that $\Bcal_{-1}$ is empty. The classes $\Gcal_d,\Bcal_d,\Scal_d$ of {\bf geometric graphs}, {\bf balls} 
and {\bf spheres} are now defined inductively as follows: $\Gcal_d$ is the set of graphs for which all unit spheres 
$S(x) \in \Scal_{d-1}$ in which case $x \in {\rm int}(G)$ or $S(x) \in \Bcal_{d-1}$ in which case $x \in \delta G$. 
The set $\Scal_{d}$ consists of graphs for which removing one vertex produces 
a graph in $\Bcal_{d-1}$. Finally, $\Bcal_{d}$ is the set of contractible graphs $G \in \Gcal_d$ for which 
the boundary $\delta G$ is in $\Scal_{d-1}$. A graph $G$ is {\bf contractible} if there exists a vertex such 
that $S(x)$ and $G-\{x\}$ are both contractible. \\

With these definitions, all graphs in $\Bcal_d$ satisfy 
$\chi(G)=1$ for $d \geq 0$ and all $d$-spheres $G$ have Euler characteristic 
$\chi(G) = 1+(-1)^d$ for dimension $d \geq -1$. On can characterize $\Bcal_d$ in $\Gcal_d$ as the 
class of graphs which admit a function with exactly one critical point and $\Scal_d$ in $\Gcal_d$.
For $d \geq 0$, $\Scal_d$ is the class of graphs in $\Gcal_d$ for which the minimal 
number of critical points is exactly $2$.

\definition{
Let $\Ccal_d$ denote the set of graphs with chromatic number $d$.
The class $\Ecal_d = \Scal_d \cap \Ccal_{d+1}$ is the class of {\bf Eulerian spheres},
and $\Dcal_d = \Bcal_d \cap \Ccal_{d+1}$ is the class of {\bf Eulerian disks}.
}

{\bf Examples.} \\
{\bf 1)} $1$-spheres are cyclic graphs $C_n$ with $n \geq $. 
Such a graph is a geometric $1$-sphere, if and only if $n$ is even. \\
{\bf 2)} A graph is a $2$-sphere if and only if it is both 4-connected and maximally planar. 
\cite{knillgraphcoloring}. The 4-connectivity implies 3-connectivity which by the Steinitz theorem 
implies that the planar graph can be realized as a convex
polyhedron in $R^3$. The maximal planarity then implies that this polyhedron has triangular 
faces. Examples are the octahedron or icosahedron. A $2$-sphere is Eulerian if an only if it
is an Eulerian graph, hence our nomenclature. \\
{\bf 3)} The $16$-cell is an example of an Eulerian $3$-sphere. 
But not all $3$-spheres are Eulerian spheres in our sense: the $600$-cell is an 
example of a $3$-sphere which is not an Eulerian sphere. Its chromatic number is $5$.  \\
{\bf 4)} Removing a vertex from an Eulerian sphere produces an Eulerian disk. 
Every interval graph is an Eulerian disk as we can color an interval with $2$ colors. \\
{\bf 5)} Gluing two Eulerian disks in $\Dcal_2$ along the circular boundary produces an 
Eulerian sphere. An Eulerian disk does not need to be an Eulerian graph. An 
example is the wheel graph $W_{6}$ with boundary $C_6$. It can be colored with $3$ colors
but the boundary points in a wheel graph have degree $3$ which is odd so that wheel graphs
are never Eulerian graphs.  \\

\begin{figure}[h]
\scalebox{0.14}{\includegraphics{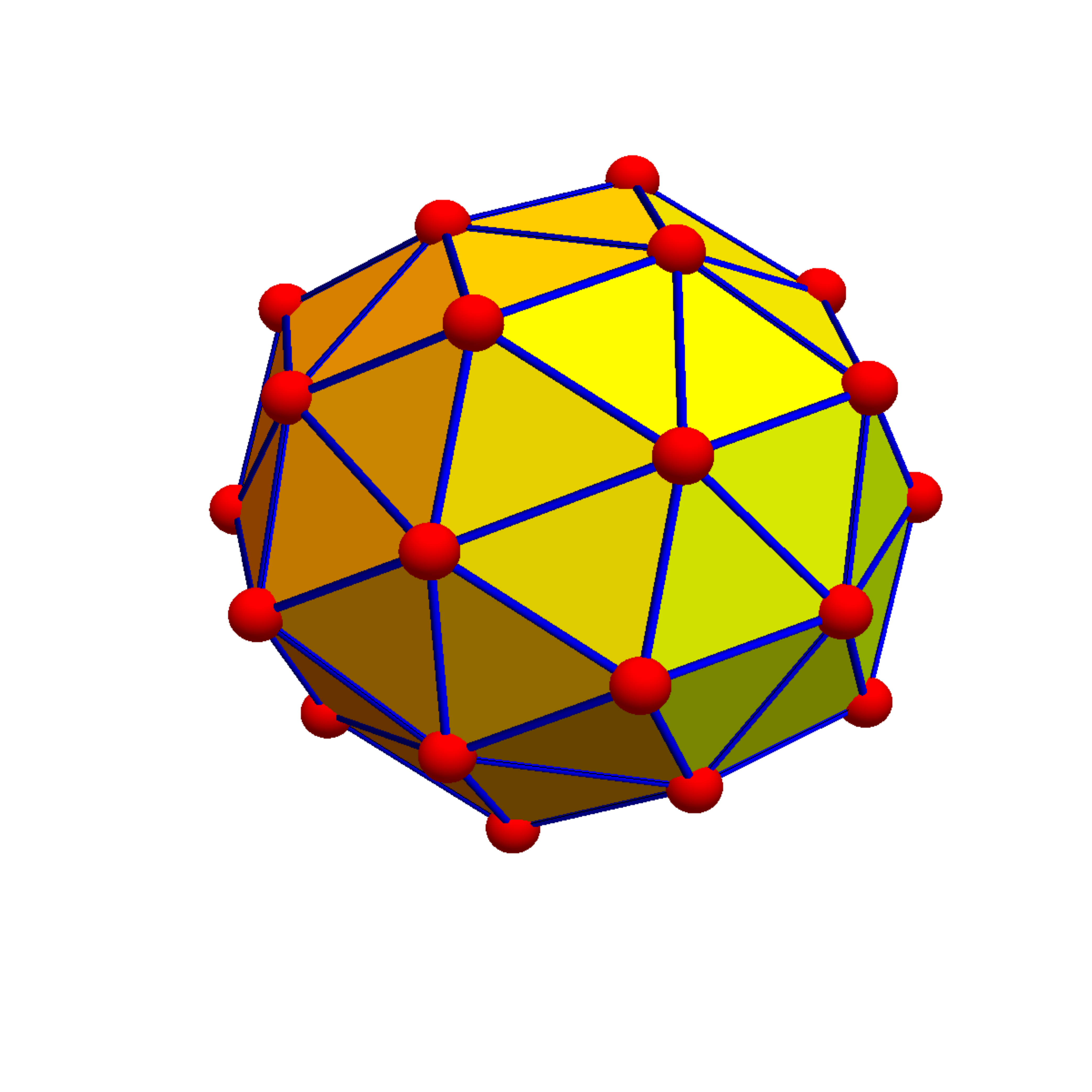}}
\scalebox{0.14}{\includegraphics{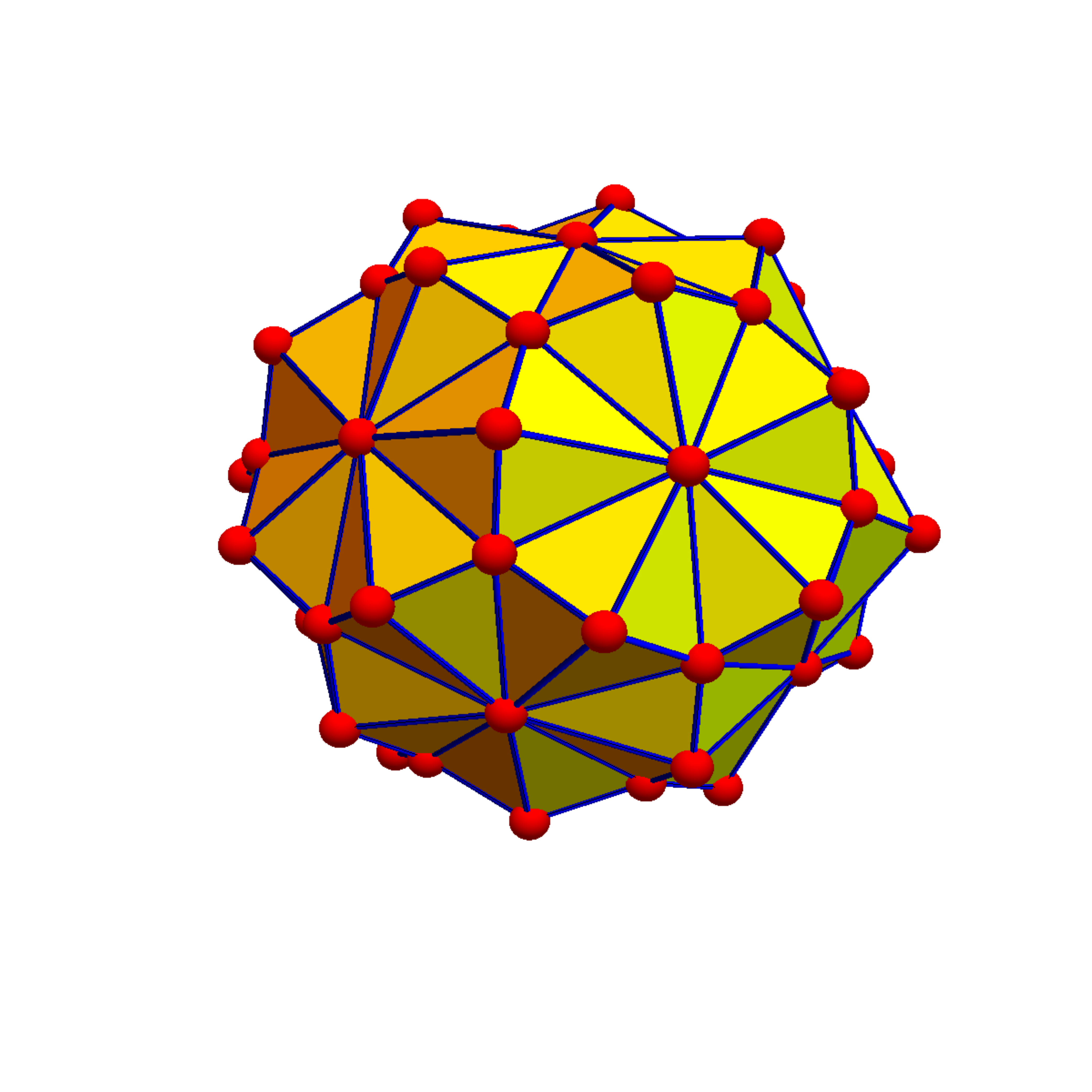}}
\scalebox{0.14}{\includegraphics{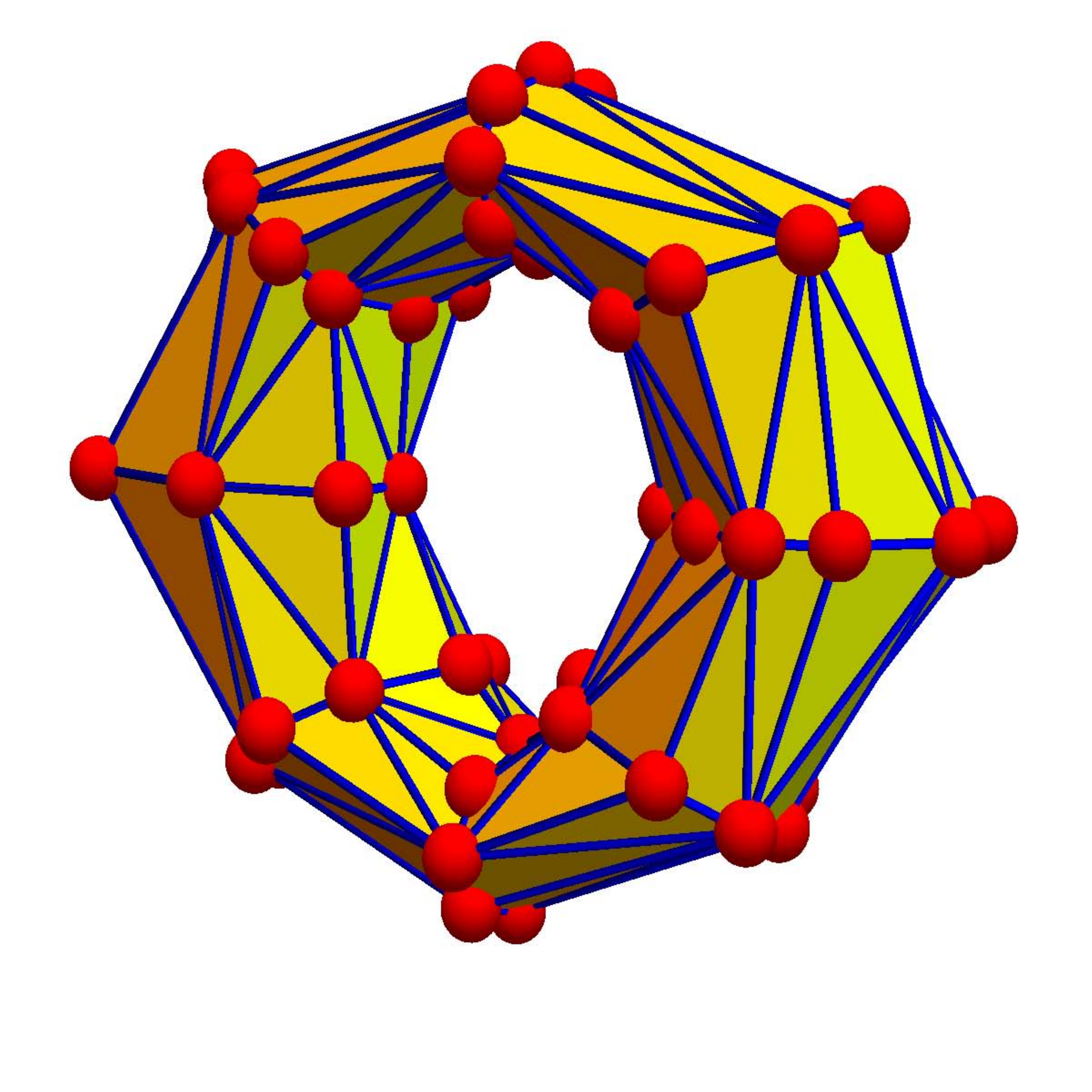}}
\scalebox{0.14}{\includegraphics{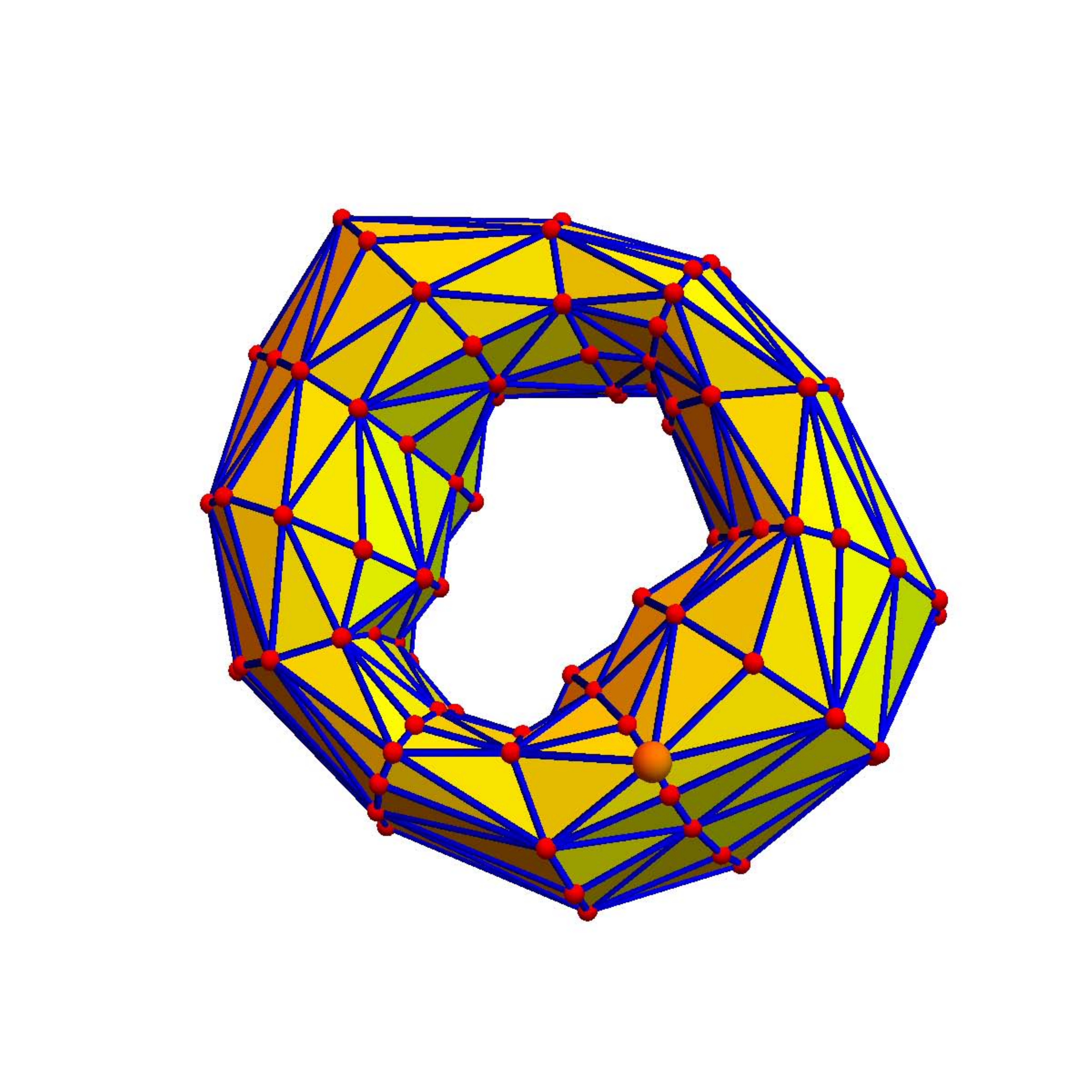}}
\caption{
The first graph is the Pentakis dodecahedron. It is a Catalan solid in $\Scal_2$ but not in $\Ecal_2$.
It is the dual of the probably the most famous semi-regular polyhedron, the {\bf "Sokker ball"} which is 
the truncated icosahedron. The second graph is a Catalan polyhedron in $\Scal_2$ which is Eulerian.
The third graph is a geometric graph in $\Gcal_2$ which is not in $\Scal_2$ but it is Eulerian.
The fourth is is an other geometric torus graph in $\Gcal_2$ which is not in $\Scal_2$ and not Eulerian
as it has 4 vertices of degree $7$. It is in $\Ccal_4$. Fisk has given a construction for
torus graphs in $\Ccal_5$. It can be obtained by surgery from this one. Just make sure that only 
two vertices of odd degree exist and that they are neighboring. }
\end{figure}

Lets start with the trivial case $d=1$ as it matches the following statements. 
For $G \in \Scal_1$, define the degree of the empty graph as the number of vertices in $G$.
We can look in general at $G$ itself as the dual graph of the empty graph, which 
can for $d=1$ be seen as a $(d-2)=(-1)$-dimensional sphere within $G$. Lets formulate
this special case too as it will matches the higher dimensional situation: 

\resultlemma{
A $1$-sphere is an Eulerian sphere if and only if every $(-1)$-dimensional subgraph 
(=the empty graph) has even degree. }

The following result is the classical Euler-Hierholzer result: 

\resultlemma{
A $2$-sphere is an Eulerian sphere if and only if every
$0$-dimensional subgraphs (=vertex) has even degree. }

{\bf Examples.} \\
{\bf 1)} The octahedron is Eulerian. \\
{\bf 2)} The stellated $2$-cube with vertex degrees $4$ and $6$ is Eulerian. \\
{\bf 3)} The icosahedron has vertex degrees $5$ and is not Eulerian.  \\
{\bf 4)} Take two Eulerian spheres which contain both a vertex $x$ or $y$ of degree $2n$.
By the connected sum construction, we can glue the two spheres along 
the unit ball to get a larger Eulerian sphere. \\

In \cite{knillgraphcoloring}, we have
defined the {\bf edge degree} of $S \in \Scal_3$ 
as the number of tetrahedra attached to an edge $e$. The union of the vertices of these tetrahedra
minus the vertices of the edge form a circular graph. 

\resultlemma{
A $3$-sphere is an Eulerian sphere if and only if every $1$-dimensional subgraphs (edge)
has even degree.}

{\bf Examples.} \\
{\bf 1)} The 16 cell, the three-dimensional analogue of the octahedron is Eulerian. 
Its clique volume vector $\vec{v}$ is $(8,24,32,16)$. \\ 
{\bf 2)} The 600 cell, the three dimensional analogue of the icosahedron is not Eulerian.  \\
{\bf 3)} The stellated $3$-cube is the {\bf tesseract}, where the $3$-dimensional 
"solid cube faces" are stellated. It is a graph containing $v_0=48$ vertices, $v_1=240$ edges, 
$v_2=384$ triangles and $v_3=192$ tetrahedra. It is an Eulerian graph for which all edge degrees
are either $4$ or $6$. \\
{\bf 4)} Using connected sum constructions, stacking together such cubes along common 
faces, we can get larger and larger classes of such graphs and get graphs in $\Gcal_3$ which 
have arbitrary large $H^2(G)=H^1(G)$ but which are minimally colorable with $d+1$ colors. \\
{\bf 5)} If $G$ is an Eulerian $2$-sphere, then the double suspension construction produces an 
Eulerian 3-sphere. For example, the $16$-cell is obtained from the octahedron by a double suspension. \\

\resultlemma{
For $d \geq 2$, the unit sphere of an Eulerian sphere is an Eulerian sphere.  }

\begin{proof}
If some $S(x)$ would need more than $d$ colors, then $G$ would need
more than $d+1$ colors to be colored. 
\end{proof} 

\definition{
The number of $k$-dimensional simplices in $G$ is denoted $v_k(G)$. 
It is called the {\bf $k$-volume} of $G$. }

The $k$-volumes $v_k$ are known to be a basis for all isomorphism invariant additive 
functionals on graphs satisfying $v_k(K)=1$ for $K=K_{k+1}$. 
The Euler characteristic is as a super volume an example of a linear combination of the volumes. 

\resultcorollary{
For an Eulerian sphere, all $k$-volumes $v_k$ are even. 
}
\begin{proof}
Use induction with respect to $d$ and use the handshaking formula
$$ v_k = (k+1)^{-1} \sum_{x \in V} V_{k-1}(x)  $$ 
which in the case $k=1$ agrees with the classical Euler handshaking result telling that
twice the number of edges is the sum of the number of edge degrees.
Since by induction, all $V_k(x)=v_k(S(x))$ are even for $k \geq 1$, the result follows
if we can show that $v_0$ is even. But since the Euler characteristic of spheres is either
$0$ or $2$ and so even, the evenness of $v_0$ follows from the formula 
$\sum_k (-1)^k v_k = \chi(G) \in \{0,2 \}$. 
\end{proof} 

One could ask whether the reverse is true and whether the condition that
all the volumes are even, assures that $G$ is Eulerian. But this is not the case.
Make two edge refinements at different but intersecting edges for $G \in \Ecal_2$. 
Each edge refinement changes the parity of the number of vertices and edges
but not for triangles. After the two refinements we have two vertices of
odd degree so that the new graph is no more Eulerian, even so all volumes are even. 
We do not believe that one can read off the Eulerian property from volumes. But one
can still ask: 

\question{
Are there some general conditions on volumes $v_k$ which force $G$ to be Eulerian?
}

Of course we do not mean silly conditions like for $d=3$, where $v_0=8$ forces the 
graph to be a 16-cell and so to be Eulerian. One could imagine a computer 
building lots of random geometric graphs and color in the grid $N^{d+1}$ 
the points $(v_0, \dots ,v_d)$ given by gometric graphs with one color, the others
with an other color. It looks like a formidable problem also to find the number
of geometric graphs in $\Gcal_d$ which have a given volume vector $(v_0, \dots, v_d)$.  

\definition{
Let $\Ecal$ denote the set of {\bf Eulerian graphs}, graphs 
on which one can find an Eulerian circuit. }

Eulerian spheres are Eulerian graphs: 

\resultcorollary{
For $d \geq 1$, we have $\Ecal_d \subset \Ecal \cap \Scal_d$.  }

\begin{proof}
By Euler-Hierholzer, we only need to know that every unit sphere has an even
number of vertices. But that follows from the previous corollary. 
\end{proof}

Because every unit sphere $S(x)$ in $G \in \Scal_d$ must be minimally colorable with $d$ colors,
the class $\Ecal_d$ of Eulerian spheres could be defined recursively as 
the class of spheres for which every unit sphere is in $\Ecal_{d-1}$ with the inductive assumption
for $d=1$, the spheres $C_{2n}$ are Eulerian.  \\

\resultlemma{
For Eulerian spheres, the dual graph is bipartite. 
}

\begin{proof}
Since spheres and especially Eulerian spheres are orientable, the coloring of each of 
the maximal simplices defines a sign of the permutation of the $d+1$ colors. This signature 
partitions the dual graph into two sets. Adjacent simplices have opposite signatures. 
\end{proof}

{\bf Examples.} \\
{\bf 1)} A $1$-sphere $C_n$ has the same dual graph and is bipartite
if and only if $n$ is even. \\
{\bf 2)} For the octahedron graph, the dual graph is the cube graph. It is bipartite. \\
{\bf 3)} There are projective planes with chromatic number $3$. In that case the 
dual graph is not bipartite. 

\resultcorollary{
The class of Eulerian spheres agrees with the class of 
spheres for which the dual graph is bipartite.}
\begin{proof}
If $G$ is an Eulerian sphere, then by the previous lemma, we have a bipartite 
dual graph. If $G$ is bipartite we can constructively color the graph with 
$d+1$ colors. 
\end{proof}

This does not generalize to non-simply connected graphs as 
the coloring will then also depend on holonomy conditions. The $6 \times 6$ torus 
for example is in $\Ccal_3$ but the $8 \times 8$ torus is not. For spheres $\Scal_d$ which 
are simply connected for $d \geq 2$, these difficulties are absent. \\

We have seen that every Eulerian sphere is an Eulerian graph and that for $d=2$,
Eulerian spheres agree with Eulerian graphs.
Of course, for $d=1$, every sphere is an Eulerian graph but only graphs of the form 
$C_{2n}$ are Eulerian graphs.  

\question{
Is there an example of a sphere for $d \geq 3$ which is an Eulerian graph but which is not
an Eulerian sphere?  }

We believe the answer is ``yes" and that $\Ecal_d \neq \Scal_d \cap \Ecal$ holds in general 
but that examples might only exist in 4 or higher dimension. We tried 
to explore the answer "no" by using an Eulerian path visiting all edges of $G \in \Scal_d$ and 
use this path to color all the vertices with a $d+1$ colors. But more likely is that there 
is a counter example. It would be surprising if $\Ecal_d = \Scal_d \cap \Ecal$ because then, we could 
focus on vertex degrees rather than $(d-2)$-dimensional simplex degrees during 
the graph modification process which could be easier. \\

Geometric graphs have spheres as unit spheres. One can look at a larger class of 
graphs, where we just ask that unit spheres are in $\Scal_g$. Lets weaken this:

\definition{
A $d$-dimensional graph is {\bf geometric exotic} if every $S(x)$ is in $\Xcal_{d-1}$. 
Let $\Xcal_d$ be the class of geometric exotic graphs. 
The induction assumption is that $\Xcal_{1}=\Gcal_1$.}

One could start even earlier and allow $\Xcal_{-1}=\{ \emptyset \}$. This would lead to an even
larger class of graphs which contains discrete varieties defined in \cite{knillgraphcoloring}.
But this would lead too far away evenso we believe that for varieties, the coloring questions are
similar to graphs in $\Gcal_d$ as singularities should not matter much. Lets go back to $\Xcal_d$
as just defined. Clearly, it follows from $\Xcal_1=\Scal_1$ that $\Xcal_2=\Gcal_2$. 
There are no exotic $2$-dimensional spaces. 
But already $\Xcal_3$ is larger as we can look for example at the double suspension of a 
higher degree surface. But this example is not homeomorphic to a standard sphere 
where homeomorphic is in the sense of \cite{KnillTopology}. 
Since we have a notion of homemorphism for graphs and because exotic spheres exist in 
classical topology, it is obvious to ask about the existence of 
``discrete exotic spheres": 

\question{
Are there examples $G$ of graphs in $\Xcal_d \setminus \Scal_d$ 
for which $G$ is homeomorphic to a graph in $\Scal_d$?
}

Such spheres could exist in higher dimensions. We would have to build a graph for
which some unit spheres are not spheres but for which there is a cover which has a
nerve graph $N$ which is graph isomorphic to a nerve graph of a standard sphere. 
If the answer to the above question is yes, one has course to wonder whether 
there is a relation with exotic spheres in topology, examples of spheres which are
homeomorphic to a unit sphere in Euclidean space but not diffeomorphic to it. 

\section{Complementary dual graphs}

The following definition is done for any subgraph $H$ of a finite simple graph $G=(V,E)$. 

\definition{
For a subgraph $H$ of a graph $G$, the intersection of all unit spheres $S(y)$ with $y \in H$
is called the {\bf complementary dual graph} of $H$ in $G$. For simplicity we write ``dual graph".}

{\bf Examples.} \\
{\bf 1)} The dual graph of the empty graph in $G$ is the graph it self. The dual graph of $G$ is the empty graph. \\
{\bf 2)} The dual graph of a $1$-vertex subgraph $x$ is the unit sphere $S(x)$ of $x$ in $G$. \\
{\bf 3)} The dual graph of a complete graph $K_k$ in $K_{d+1}$ is the complete graph $K_{d+1-k}$
         formed by the complementary vertices of $K_k$. \\

\begin{figure}[h]
\scalebox{0.13}{\includegraphics{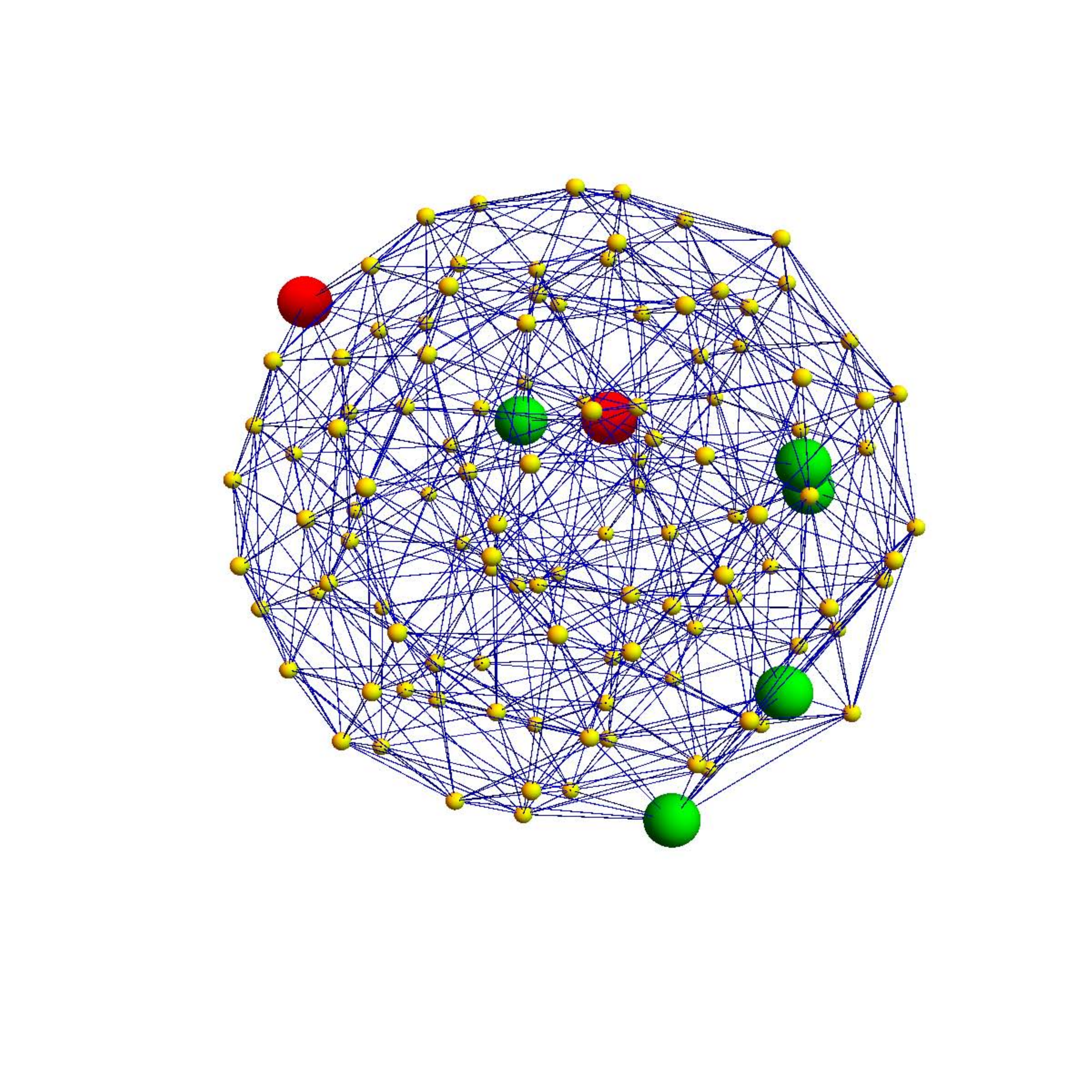}}
\scalebox{0.13}{\includegraphics{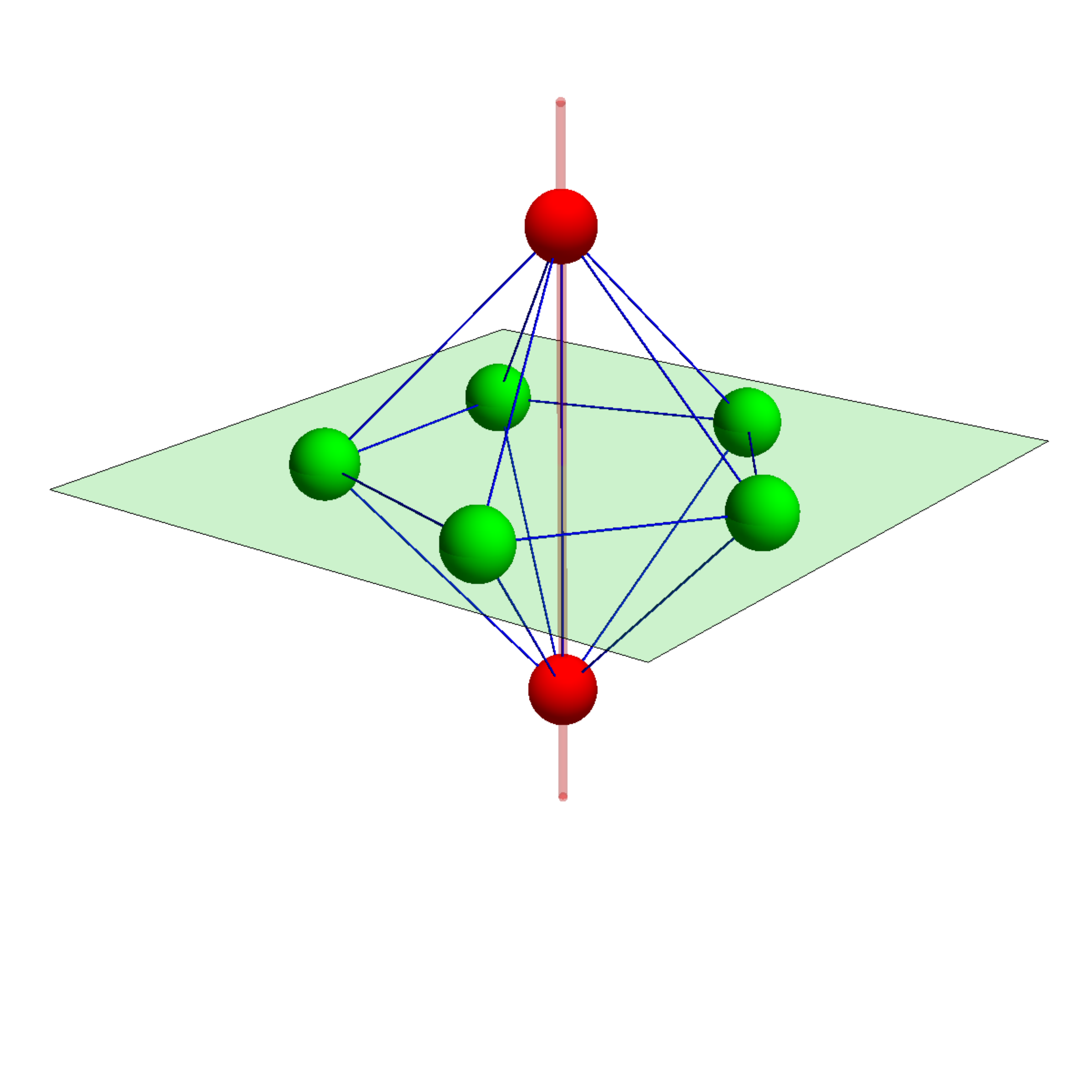}}
\caption{
A graph with two vertices and its dual graph. The graph is the $600$-cell, a $3$-sphere.
The geometric picture to the right shows two spheres. One is the 
intersection of a line with the unit sphere and the second is the 
intersection of a plane with the unit sphere. }
\end{figure}

\begin{figure}[h]
\scalebox{0.13}{\includegraphics{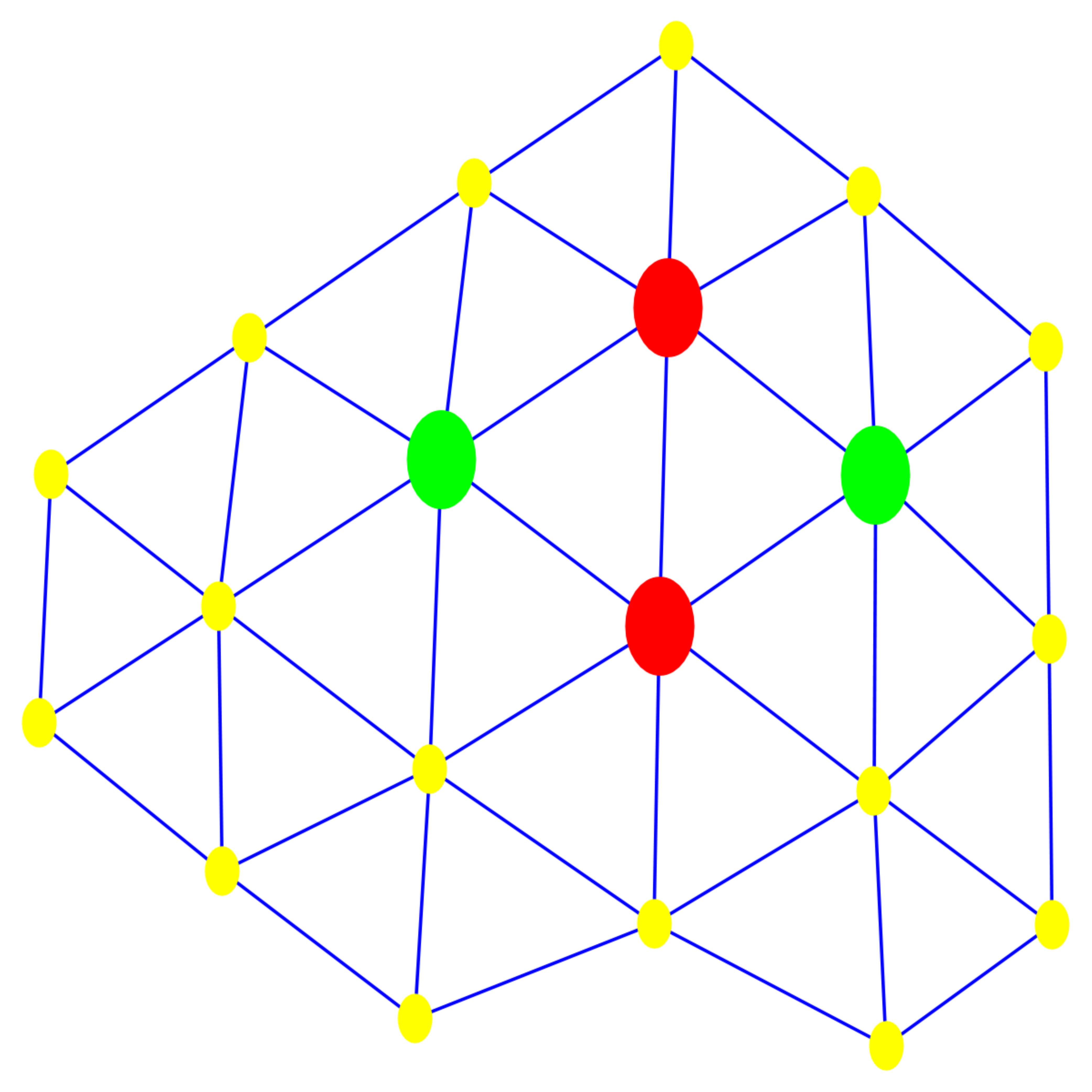}}
\scalebox{0.13}{\includegraphics{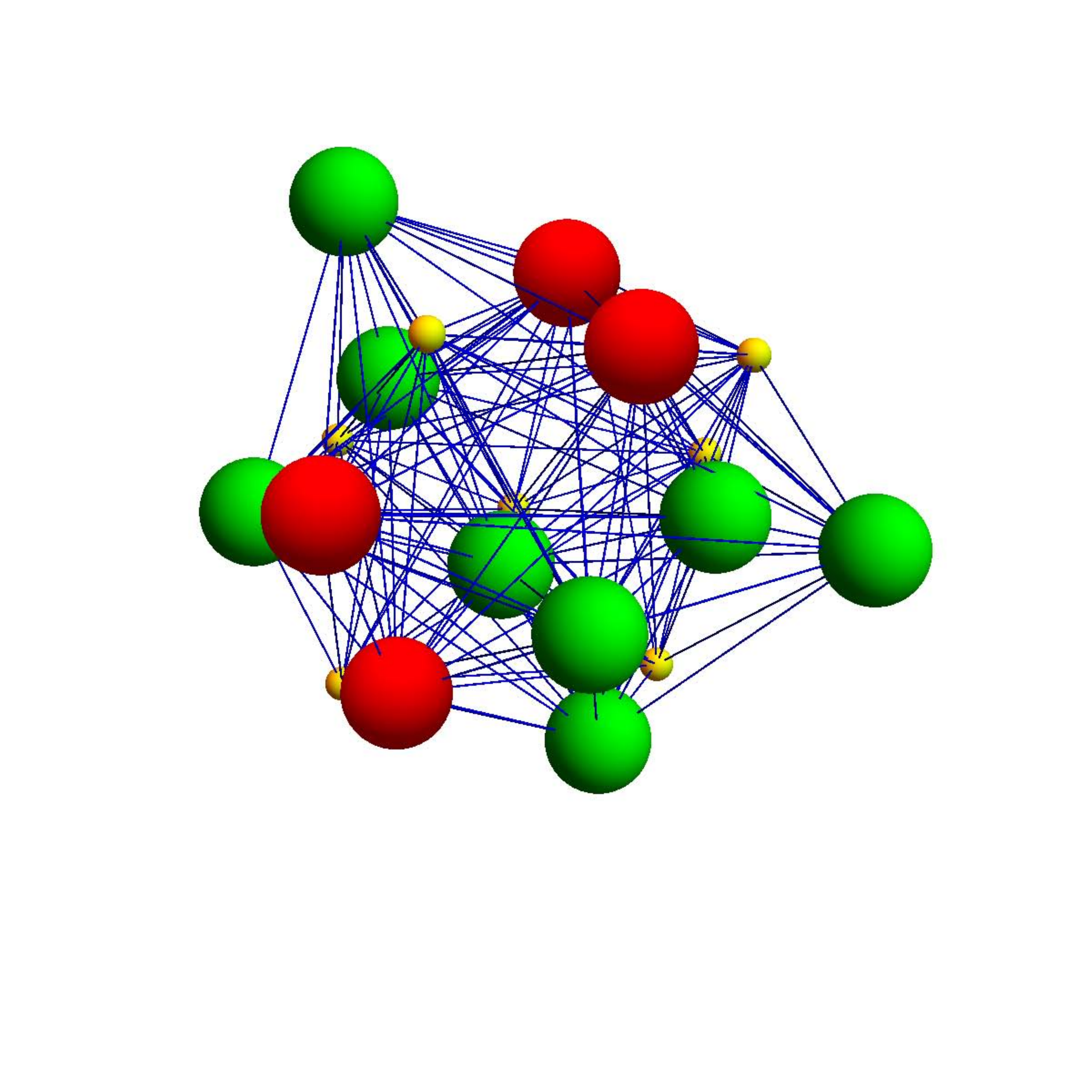}}
\caption{
A graph $H=K_2$ in a hexagonal graph $G$ and its dual graph $\hat{H}$
which is a sphere in $\Scal_0$. The right figure shows a random Erdoes-R\'enyi graph 
and the dual graph of a subgraph $H$  with $3$ elements. The graph $H$ has been obtained by first
computing the dual graph of an other graph. According to the lemma, we are locked in
then to a duality. 
}
\end{figure}

\resultlemma{
For a complete subgraph $K_k$ of a $d$-dimensional geometric graph, the 
dual graph is a $(d-k)$-sphere. }

\begin{proof}
For $k=1$, the complete graph $K_1$ is a single vertex and the dual graph is its unit sphere, 
which is a $(d-1)$-sphere: the empty graph. In general, use induction with respect to to $k$: given a subgraph
$K_{k+1}$, take a point $x$ off from $K_{k+1}$ to get a complete subgraph $K_k$. 
The dual graph of $K_k$ is now a $(d-k)$-dimensional graph $S$ containing $x$. 
By induction assumption, it is a sphere. 
The dual graph of $K_{k+1}$ is the unit sphere of $x$ in $S$ which (because $S$ is a sphere)
by definition is again a $d-k-1 = (d-(k+1))$-dimensional sphere. 
\end{proof} 

{\bf Examples.} \\
{\bf 1)} For a vertex $x$ in the octahedron, the dual graph $\hat{x}$ is the $1$-dimensional
unit sphere $S(x)$ which is $C_4$ in this case. \\
{\bf 2)} For an edge $e=(a,b)$ in a $2$-dimensional geometric graph, the dual graph $\hat{e}$ 
consists of the zero-dimensional sphere consisting of two edges $\{ c,d \; \}$, where 
$(a,b,c)$ and $(a,b,d)$ are the two triangles adjacent to $e$. \\

The following statement justifies the name duality: 

\resultlemma{
Let $G$ be a general finite simple graph. 
For any subgraph $H$ of $G$, one has $H \subset \hat{\hat{H}}$. 
If $H$ is a dual graph in $G$, then $\hat{\hat{H}} = H$. }

\begin{proof}
Given a vertex $x \in H$. Every vertex $y$ in $\hat{H}$ has distance $1$ from $x$. 
Therefore $x$ is in the dual graph of $\hat{H}$. This shows that 
$H \subset \hat{\hat{H}}$. 
Applying this to $\hat{H}$ shows that $\hat{H} \subset \hat \hat \hat{H}$. (1)
We know that $A \subset B$, then $\hat{B} \subset \hat{A}$ as $\hat{A}$ consists of
less intersections. Applying this to $A = H$ and $B = \hat{\hat{H}}$ shows
that $\hat{\hat{\hat{H}}} \subset \hat{H}$ (2). 
From (1) and (2) we see $\hat \hat{H} = H$ if $H$ is a dual graph. 
\end{proof}

Is every subgraph $H$ of a graph $G$ a dual graph of some graph? 
No. Otherwise, the dual operation would be bijective but there are many 
subgraphs $H$ of $G$ for which $\hat{H}$ is the empty graph. Examples are graphs of 
diameter larger than $2$. In a geometric situation we can characterize dual graphs
as follows: 

\resultlemma{ 
For $G \in \Scal_d$, the dual graph of any subgraph $H$ is either empty, a sphere,
a complete graph, a ball or the entire graph.}

\begin{proof}
For the empty graph $H$, the dual graph $\hat{H}$ is the entire graph $G$. 
Graphs of diameter $1$ are complete subgraphs $H$ for which the dual is a sphere.
If the diameter is $2$ (within $G$) we deal with subsets of spheres $S(x)$
and the dual is a sphere together with the center which is a disk or a complete 
subgraph. For a graph of diameter $3$ or more, two spheres do not intersect and 
the dual graph is empty.
\end{proof} 

{\bf Examples.} \\
{\bf 1)} If the host graph $G$ is a complete graph, every subgraph $K$ has the 
dual graph graph $\hat{K} = G \setminus K$ and these two complete graphs are dual 
to each other. \\ 
{\bf 2)} In $G=C_4$, the graph $G$, the empty graph $\emptyset$, the
one point graphs and their spheres are the only dual graphs. The dual 
graph of a unit ball is the empty graph. \\
{\bf 3)} For a 3-sphere, the dual graph of an edge is a circular graph $S$. 
If we take a sufficiently large part $H$ of that circular graph, then the dual $\hat{H}$ 
is already the edge again. But $\hat{\hat{H}} = S$. 
For a complete subgraph $H$ it is a sphere $S$. Now complete this sphere to a ball $B(x)$.
The complementary graph of this ball is a $(d-1-k)$-dimensional                  
sphere $\hat{S}$. The two spheres $S$ and $\hat{S}$ are dual spheres in          
$S(x)$ in the sense that the dual sphere of the ball of $S(x)$ is $\hat{S}(x)$
and the dual sphere of the ball of $\hat{S}(x)$ is $S(x)$.  \\
{\bf 4)} For an edge $e=(a,b)$ in a $4$-sphere, the dual graph $\hat{e}$ is a $2$-sphere $S$, 
the intersection of two $3$-spheres $S(a)$ and $S(b)$. The dual of $S$ is again $e$. \\
{\bf 5)} For a vertex $x$ in $\Scal_d$  the dual graph $\hat{x}$ is a sphere 
and the dual of $\hat{x}$ is again $\{x\}$ a complete graph. 

\section{Degree condition for Eulerian} 

We have seen that the degree of the empty graph $\emptyset$ is the order of the graph, 
that the degree of a vertex in a $2$-dimensional graph is 
the degree of the vertex and that the degree of an edge in a 
$3$-dimensional graph is the edge degree, the order of the dual graph of $e$. 
In general, we have:

\definition{
Given $G \in \Gcal_d$, let $x$ be a $k=(d-2)$-dimensional simplex in $G$. The 
{\bf degree} of $x$ is defined as the order of the $1$-sphere $\hat{x}$. }

Here is the main result. It generalizes what we know already for $d=1,d=2$ and $d=3$:

\resulttheorem{
A sphere $G \in \Scal_d$ is an Eulerian sphere if and only if 
every $(d-2)$-dimensional simplex $x$ in $G$ has even degree.}

\begin{proof}
Let $x$ be a $(d-2)$-dimensional simplex and let $\hat{x}$ be its dual sphere
which we know to be a $1$-dimensional circular graph. If $n$ is the number of 
vertices in $\hat{x}$, it is the degree of $x$. 
Let $y$ be a $d$-dimensional simplex which contains $x$. 
It is given by the vertices in $x$ as well as two adjacent vertices $(a,b)$
in $\hat{x}$. It follows that the set of $d$-dimensional simplices which hinge
at $x$ form a circular chain of length $n$. Now pick a simplex $z$ in this
set and color it with $d+1$ colors. The coloring of the adjacent simplices
in the circular chain is now determined. Since $n$ is even we can continue through 
and color all the $n$ simplices with $d+1$ colors. Now pick an other 
$d-1$ dimensional simplex which is the intersection of two simplices. It again
produces a chain of simplices. In the same way, we can continue the coloring. 
Since the graph $G$ is simply connected, we can color the entire graph with 
$(d+1)$ colors. 
\end{proof}

\begin{figure}[h]
\scalebox{0.09}{\includegraphics{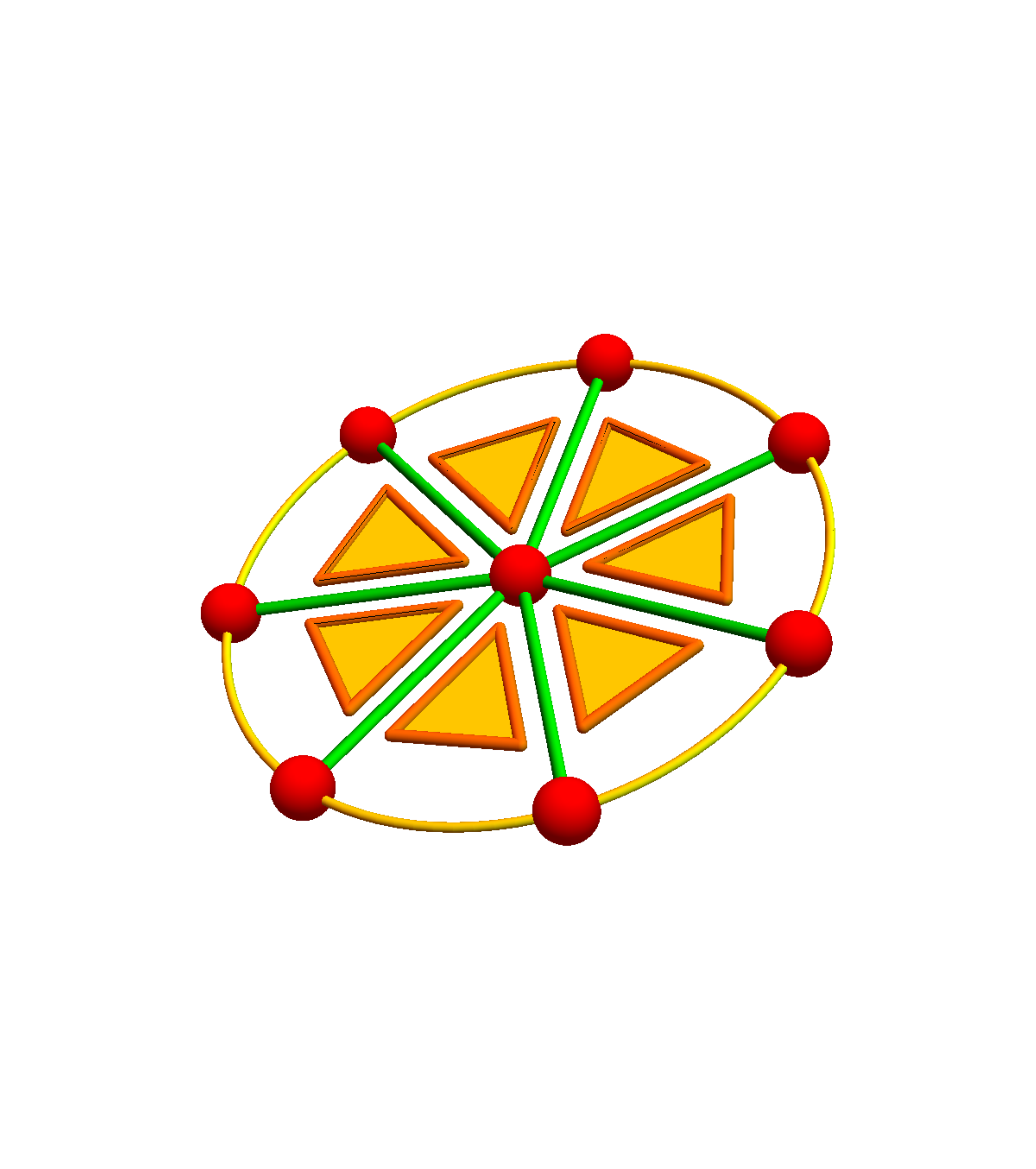}}
\scalebox{0.09}{\includegraphics{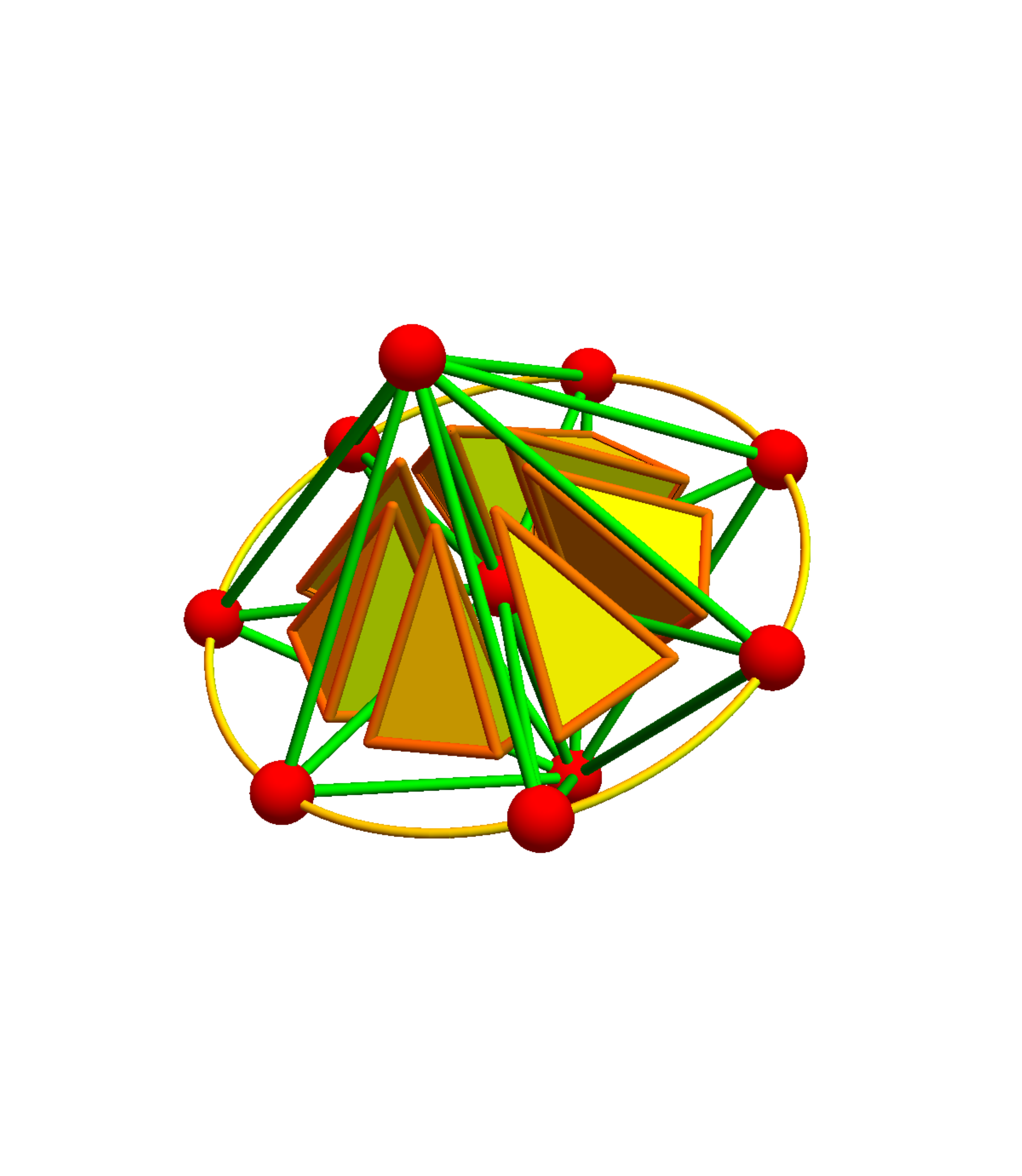}}
\scalebox{0.09}{\includegraphics{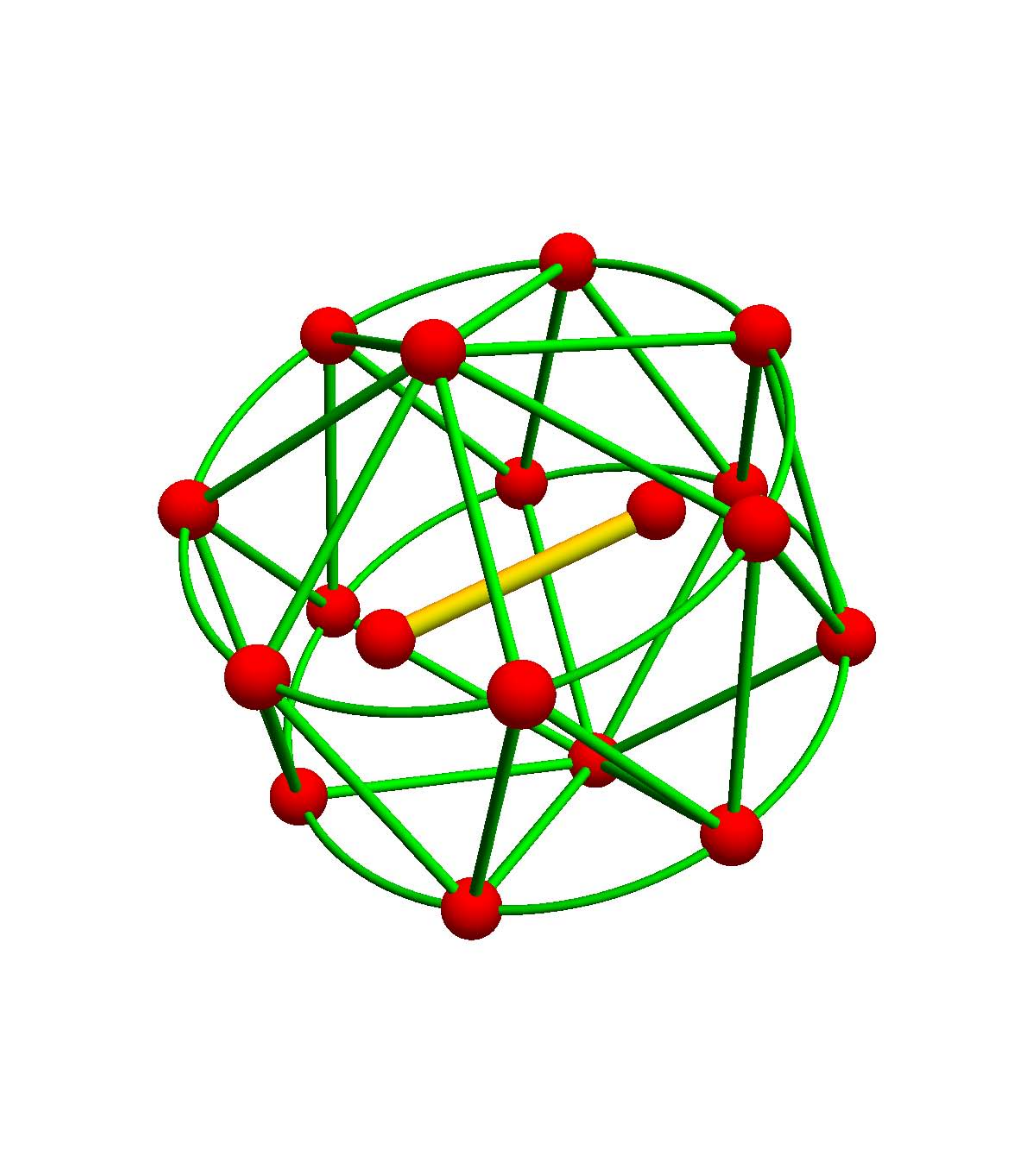}}
\caption{
The degree of a vertex in a sphere $G \in \Scal_2$ is
the degree of the vertex. The degree of an edge in a sphere $G \in \Scal_3$
is the edge degree. The dual graph of an edge $e$ in a sphere $G \in \Scal_4$
is a 2-sphere. The degree of each triangle in this graph is defined as the
number of hyper-tetrahedra $K_5$ attached to the triangle. }
\end{figure}

\section{Edge refinement and collapse}

There are many refinements of Eulerian graphs which 
keep the graph Eulerian. The already mentioned connected sum construction 
can be seen as such: chose a unit ball $B(x)$ inside $G$ and glue two
copies of the graph together along the boundary. This can be seen as a refinement,
as we can see the second graph glued into the unit ball of the first. We are more
interested in refinements which do modify the Eulerian structure as we want
to get to the Eulerian graphs eventually. One example is to do an edge division and 
the inverse, the edge collapse. \\

Given a $d$-sphere $G$ and a $(d-2)$-dimensional complete subgraph $H$.
The dual graph $\hat{H}$ is a cyclic graph. We want to find a subdivision
algorithm which renders all degrees even. 

\definition{
An {\bf edge subdivision step} for $G \in \Scal_d$ consists of taking
an edge $e=(a,b)$, subdividing it with a new vertex $x$ and
connecting $x$ with all the vertices in the dual graph $\hat{e}$. }

The reverse of an edge division can be identified with edge collapse,
which is often used in graph theory but which does not always preserve
dimension and so the class of geometric class. 

\definition{
An {\bf edge collapse} with an edge $e=(a,b)$ identifies the two vertices $a,b$.
If $x$ is the newly identified point, all edges $(y,a)$ become $(y,x)$ and
all edges $(y,b)$ become $(y,x)$.}

Each edge subdivision preserves the class of spheres and is an Ivashchenko 
homotopy:

\resultlemma{
If $H$ is obtained from $G \in \Scal_d$ by an edge $e$ subdivision step, then 
$H$ is again in $\Scal_d$ and is homotopic to $G$. The graph $G$ can be obtained
back from $H$ by an edge collapse. }

\begin{proof}
(a) {\bf $H \in \Gcal_d$}:  \\
This is seen by induction with respect to $d$. 
The sphere $S(x)$ of the new vertex $x$ dividing the edge $e=(a,b)$ is a
double suspension of the sphere $S(a) \cap S(b)$ and so a sphere. 
The spheres $S(a)$ and $S(b)$ do not change topologically. For the sphere
$S(a)$, the vertex $b$ is just replaced with the vertex $x$. 
The other affected vertices are vertices in $S(a) \cap S(b)$. 
For such a vertex $y$, the sphere contains originally the edge $(a,b)$
and afterwards the edge $(a,x)$ and $(b,x)$ where $x$ is replaced with
any vertex in $S(y) \cap S(a) \cap S(b)$. In other words, the sphere $S(y)$
has undergone an edge refinement too and by induction this is again is
a sphere.  \\
(b) {\bf $H$ is homotopic to $G$.} \\
The edge division can be split up into two parts. First by
adding a new vertex $x$ and attaching it to the contractible 
subgraph $S(a) \cap S(b) \cup \{a,b\}$.
This does not leave spheres invariant but it is a homotopy step. 
Now remove the edge $(a,b)$. This is possible since $S(a) \cap S(b)$
is now contractible. By the way, this is the original Ivashchenko notion of 
I-contractibility which additionally included that deformation. \\
(c) {\bf $H \in \Scal_d$}: since by b), the graph is homotopic to a sphere, it must be
in $\Scal_d$.  \\
(d) {\bf Edge collapse}. If $e=(a,b)$ and $x$ was the newly added vertex
in the edge, then either collapsing the new edge $(x,a)$ or collapsing $(x,b)$
produces a graph isomorphic to the old graph $G$. 
\end{proof} 

Note that unlike edge refinements and the fact that edge collapse
reverses an edge refinement, an edge collapse in general is not a homotopy. 
An example is the graph $C_4$ which is not simply connected and has Euler characteristic
$0$. After an edge collapse we have a triangle which is contractible and has 
Euler characteristic $1$. If we do an edge subdivision for the triangle,
we obtained the kite graph.  \\

\begin{figure}[h]
\scalebox{0.13}{\includegraphics{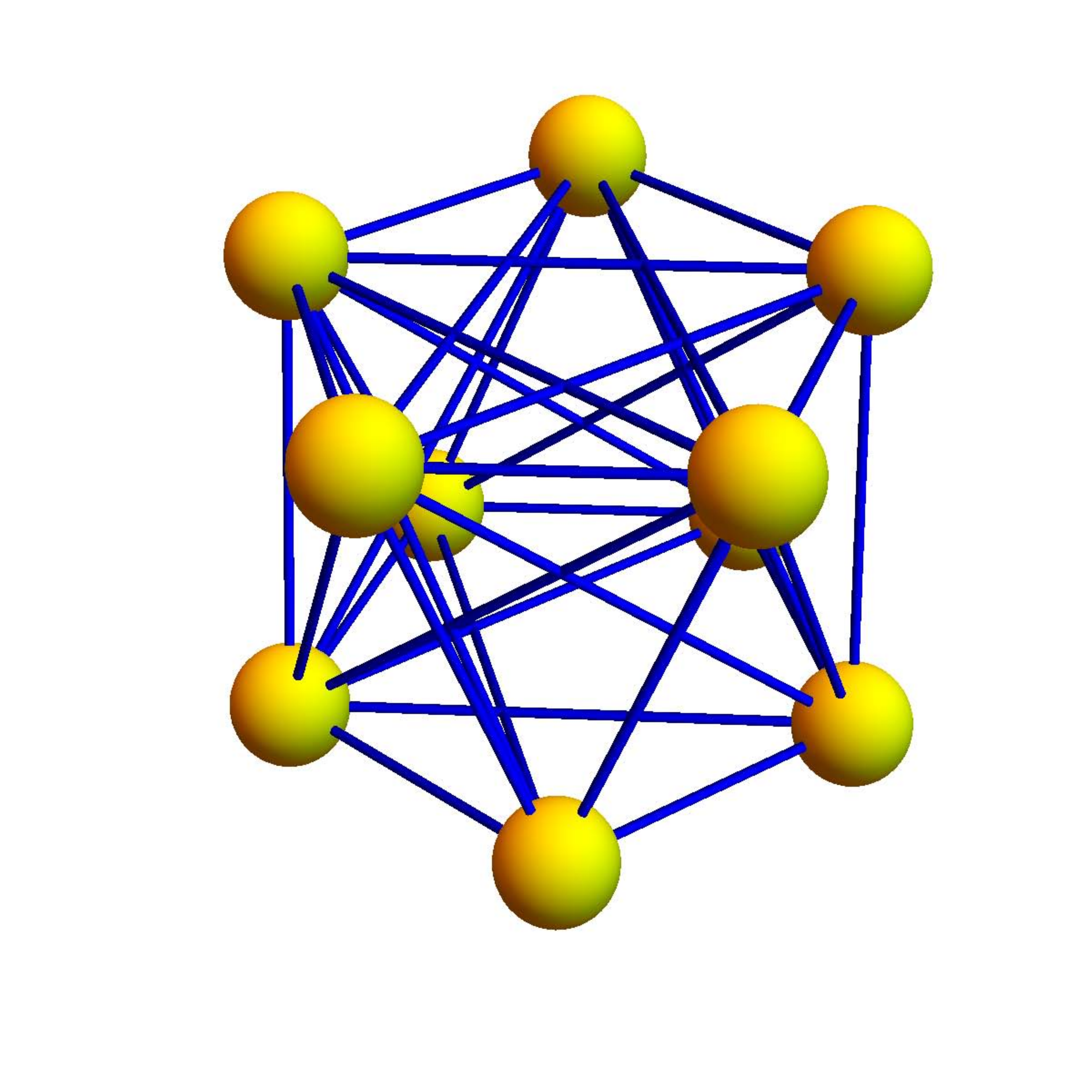}}
\scalebox{0.13}{\includegraphics{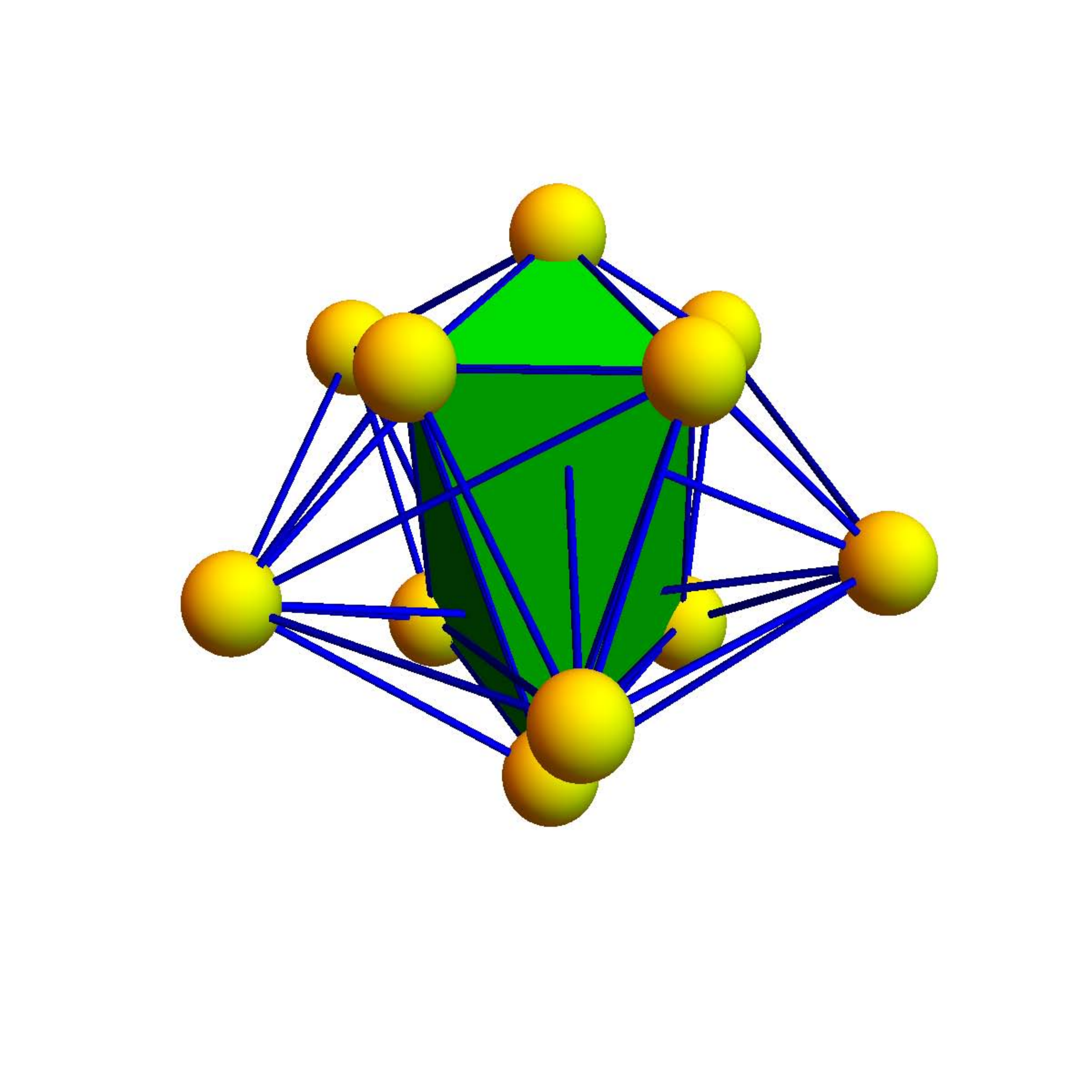}}
\caption{
A $4$-sphere $G$ with 32 four-dimensional chambers $K_5$. It is Eulerian and the
only Platonic sphere as we will see below. The clique data are $\vec{v} = (10, 40, 80, 80, 32)$.
The second picture shows the graph after an edge $e$ subdivision process has been applied.
We have shaded the triangles which emerged with odd degree. 
The triangles on a $2$-sphere $\hat{e}$ have now odd degrees.
There is 1 vertex, 7 edges, 8 triangles, 20 tetrahedra and 18 chambers more.
The deformed graph is homotopic to the first. }
\end{figure}

\begin{figure}[h]
\scalebox{0.13}{\includegraphics{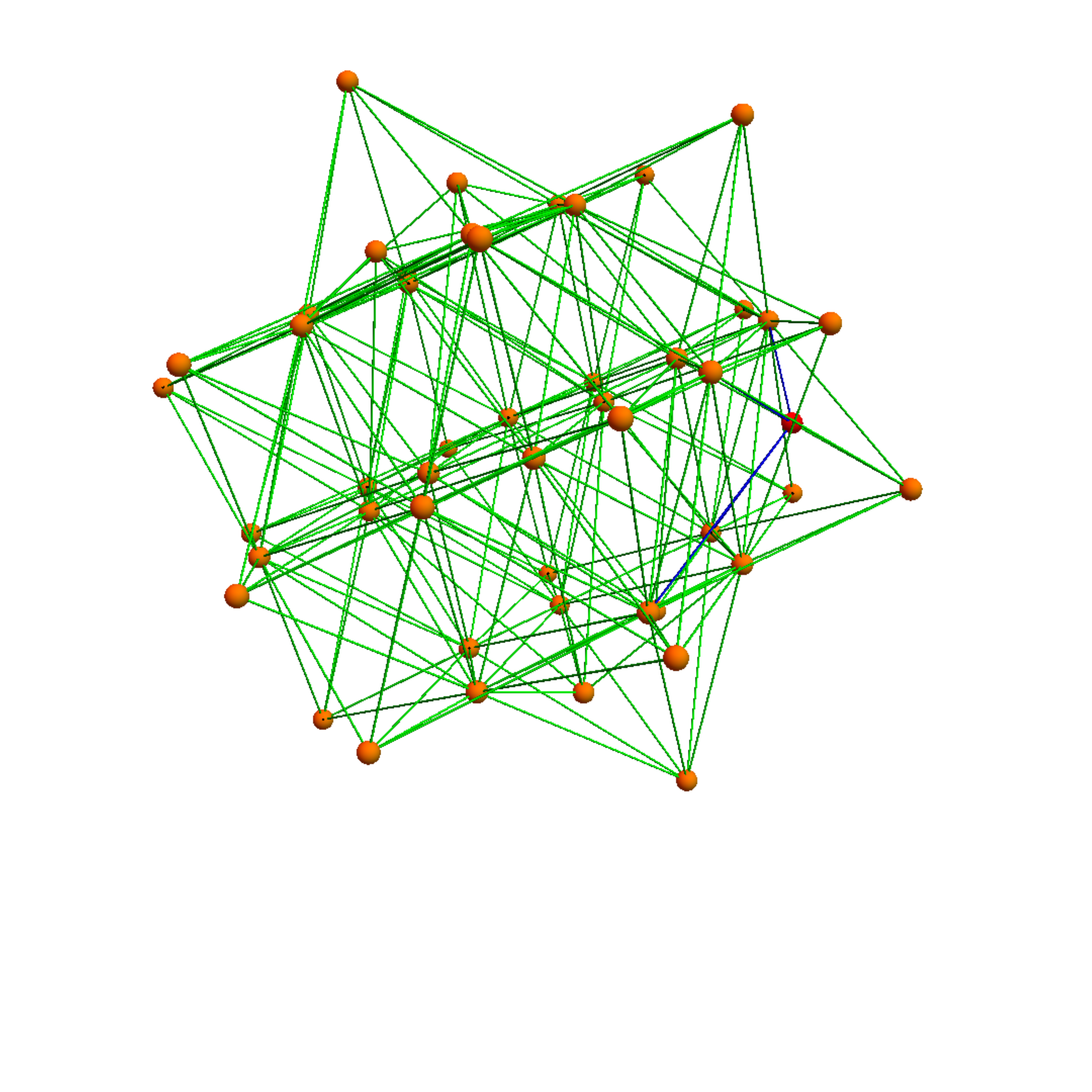}}
\scalebox{0.13}{\includegraphics{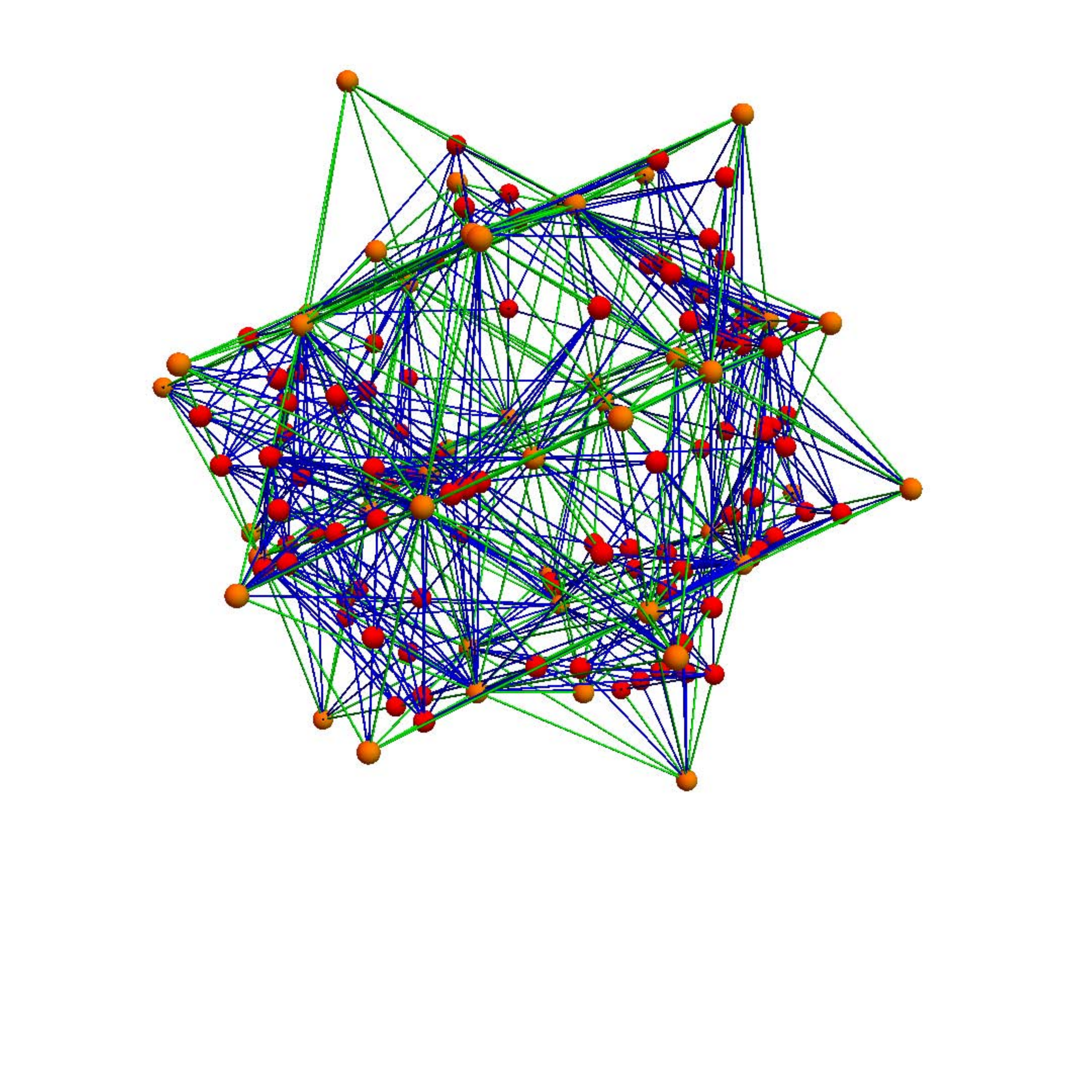}}
\caption{
A triangularization of the tesseract to the left. It is a graph in $\Ecal_3$. 
To the right we see the graph after a few random edge refinements have been applied. 
The graph to the right is still a three dimensional sphere but no more in $\Ecal_3$. 
}
\end{figure}

What is the effect of an edge division or edge collapse? 

\resultlemma{
An edge division step with edge $e$ 
increases by $1$ the degree for all $(d-2)$-simplices in $\hat{e}$. }

\begin{proof}
Given a $(d-2)$ simplex $z$ in $\hat{e}$, it gives together with the edges
$(a,b)$ a $d$-dimensional simplex which counts to the degree of $z$. 
When splitting $e$ to $(a,x,b)$, we get two simplices, one with $(a,x)$
added to the vertices of $z$, the other with $(x,b)$ added. 
\end{proof}

{\bf Examples.} \\
{\bf 1)} For $d=1$, the dual graph $\hat{e}$ of any edge $e$ is empty. 
There are no $(d-2)$-simplices in the graph.  
The subdivision changes the degree of the entire graph. This is a
special case; not at least because $1$-dimensional spheres are not simply connected. \\
{\bf 2)} For $d=2$, the dual graph $\hat{e}$ consists of two vertices. The 
subdivision changes the degree of these two vertices. \\
{\bf 3)} For $d=3$, the dual graph $\hat{e}$ consists of a circular graph. 
The subdivision changes the degree of each of the edges of this circle. \\
{\bf 4)} For $d=4$, the dual graph $\hat{e}$ consists of a $2$-sphere, 
a polyhedron. The subdivision changes the degrees of each triangle in that sphere. 
Remember that the degree of the triangle is the number of elements in the $1$-dimensional
circle dual to that triangle. \\

We believe that every $d$-sphere can be modified by refinements or collapses within the class $\Scal_d$ 
so that it becomes Eulerian but we only know this in the case $d \leq 2$. 

\resultlemma{ A $2$-sphere can be modified by edge refinements to become Eulerian.  }

\begin{proof}
Chose two vertices of odd degree. Find a path of even length connecting the two
and cut the triangles. Possibly first make a subdivision first. 
\end{proof}

In the following, an {\bf edge collapse} only refers to edge collapses which 
keeps the graph in $\Scal_d$. Irreducible graphs to be defined below would 
not allow an edge collapse. 
The following statement would provide a proof of the 4-color theorem:

\conjecture{
A $3$-sphere can be modified by edge refinements or collapses
to become Eulerian without modifying a given embedded 2-sphere inside.}

The following statement would provide a proof of the conjecture that every graph in 
$\Gcal_2$ can be colored with $3,4$ or $5$ colors.

\conjecture{
A $4$-ball $B$ can be modified by edge refinements or collapses  to become Eulerian without
modifying its boundary $3$-sphere $S=\delta(B)$ and an embedded two-dimensional surface $G$
placed in $B$ so that $G^c \cap {\rm int}(B)$ is simply connected.}

The following statement would provide a proof of the conjecture  \\
$\Scal_d \subset \Ccal_{d+1} \cup \Ccal_{d+2}$: 

\conjecture{
A $(d+1)$-sphere can be modified by edge refinements or collapses to become Eulerian without
modifying an embedded $d$-sphere inside. }

\section{Projective spheres and irreducibility}

Spheres in spheres play the role of linear subspaces of the tangent space in 
classical geometry: the reason is that if we take a linear subspace and intersect
it with the sphere, we get a lower dimensional sphere. Since in the discrete, we only 
have spheres and no notion of linear subspaces, we work with the later
similarly as classical geometry can be dealt with by compactifying Euclidean space
leading to projective geometry. Since we use edge subdivision to modify graphs aiming
to make them Eulerian, we are interested in graphs which can not be obtained from 
edge subdivision steps. 

\definition{
A graph $G \in \Gcal_d$ is called {\bf irreducible} if no edge collapse
can be applied to it without getting out of $\Gcal_d$.}

There is a different notion of irreducibility in discussions about the 4 color theorem but 
should be no confusion. \\

{\bf Examples.} \\
{\bf 1)} An octahedron is irreducible because it is the smallest graph in $\Scal_2$. \\
{\bf 2)} The octahedron is a double suspension of $C_4$. Any double suspension $G_n$ of $C_n$ with $n \geq 4$
is reducible. The reason is that $G_{n+1}$ is obtained from $G_n$ by an edge refinement. \\
{\bf 3)} The icosahedron is reducible as we can collapse one edge and still keep it geometric.  \\
{\bf 4)} All higher dimensional cross polytopes in $\Scal_d$ are irreducible. \\
{\bf 5)} There is exactly one irreducible sphere in $\Scal_1$: it is $C_4$. \\
{\bf 6)} Irreducibility is also defined for graphs in $\Gcal_d$. A hexagonal flat torus for example
is irreducible since an edge refinement steps produces a vertex of degree $4$.  \\

Can we classify irreducible spheres? Yes, in dimension $d=1$, where only $C_4$ is irreducible and also in 
dimension $d=2$: after an edge refinement step, we always have a vertex of degree $4$ so that if there 
are no degree $4$ vertices, we can not collapse: 

\resultlemma{
If $G \in \Scal_2$ has a vertex $x$ of degree $4$ for which the disc $\{ y \; | d(x,y) \leq 2 \}$
is in $\Bcal_2$, then we can apply an edge collapse. }
\begin{proof}
Let $a,b,c,d$ be the neighbors of $x$ and $(a,b)$ a pair of elements in $S(x)$.
Remove the edges $(x,a),(x,b),(x,c),(x,d)$ and add the edge $(a,b)$. By the radius of
injectivity condition, this edge collapse still produces a graph in $\Scal_2$.
\end{proof}

Graph modifications are important in the four color theorem. 
One can restrict to graphs in $\Scal_4$ for which all vertex degrees are
larger or equal to $4$. Kempe used his chains to avoid degree 4 vertices and overlooked
a case when trying to avoid degree-5 vertices which by Gauss-Bonnet $\sum_{x} (1-{\rm deg}(x)/6)=2$
always exist for $G \in \Scal_2$ if degree-4 vertices have been excluded. 

\definition{
A graph $G$ is called {\bf projective} if every unit sphere $S(x)$ admits a fixed point free 
involution $T(x)$. It is {\bf strongly projective} if
every $S(x)/T(x)$ is a discrete projective space. It is called {\bf weakly projective} if 
every unit sphere admits an involution $T$ for which a $d-1$ sphere is fixed. 
}

{\bf Examples.} \\
{\bf 1)} For a $16$ cell, every unit sphere is an octahedron $G$ which is projective but not strongly projective. 
It admits an involution but the quotient $G/T$ is no more geometric. There are no strictly 
positive curvature projective planes in the discrete. \\
{\bf 2)} If every cyclic graph is $G=C_{2n}, n \geq 4$, then $G$ is strongly projective. It has an involution $T$
such that $G/T$ is again a cyclic graph. \\
{\bf 3)} Graphs containing an unit sphere with an odd number of vertices can not be projective. Therefore, the degree has to be even
for all vertices. \\
{\bf 4)} Strongly projective spheres can be constructed by taking a discrete projective space $G$, for example by taking
a nice triangulation of a classical projective space, then take a double cover. \\
{\bf 5)} A decahedron, the polyhedron obtained by doing a double suspension on $C_5$ is an example
of a weakly projective space. There is an antipodal map, which assigns to each vertex $p$ the unique point
with largest distance from $p$. But this map can not be fixed point free as the restriction to the
equator $C_5$ shows, where the odd number of vertices prevent a fixed point free involution. 

\resultlemma{
The class of projective $d$-spheres is contained in the class Eulerian of $d$-spheres. 
}

{\bf Examples.} \\
{\bf 1)} A graph $G$ in $\Scal_2$ which is projective then it is in $\Ecal_2$ because all unit spheres have then even degree. \\
{\bf 2)} A graph $G$ in $\Scal_3$ is projective then every unit sphere $S(x)$ in $G$ is projective and so Euclidean. 
This means that $\Ecal_3$ contains the projective $3$-spheres. \\
{\bf 3)} There are Eulerian graphs which are not projective. A simple example is obtained by taking an octahedron and
subdividing the same edge twice. This produces a new Eulerian graph but it does not admit a fixed point free involution. 

\definition{A $d$-sphere $S$ is {\bf generic} it is not weakly projective and all its unit spheres are generic. 
A graph $G \in \Gcal$ different from a sphere is generic if all unit spheres $S(x)$ are generic. }

{\bf Examples.} \\
{\bf 1)} The icosahedron $G$ is generic in $\Scal_2$. \\
{\bf 2)} The $600$ cell $G$ is generic in $\Scal_3$. \\
{\bf 3)} Any flat torus $G$ with hexagonal tiling is generic in $\Gcal_2$. \\

In other words, $G$ is weakly projective if all unit spheres $S(x)$ are a double cover of a disc ramified over 
its boundary or even a double cover over a projective space. Two weakly projective spheres can be added 
with the direct sum construction so that $G + G = G$.  A generic sphere does not allow a nontrivial involution 
on each of its spheres. Strongly projective spheres are spheres which are double covers of projective spaces.
Any large enough projective sphere is strongly projective. 

\resultlemma{ Generic spheres are irreducible. }
\begin{proof}
We prove that if $G$ is not irreducible, then the graph $G$ is not generic. 
If $G$ is not irreducible, it is obtained from an other graph $H$ by an edge division.
This edge refinement step has added a vertex $x$ in the middle of an edge $e=(a,b)$ 
and produced a projective sphere $S(x)$ within $G$: the involution switches $x \leftrightarrow y$
and leaves the sphere $S(a) \cap S(b)$ invariant. The graph is not generic. 
\end{proof}

\resultlemma{
If $G$ is a generic sphere and $K$ is a simplex, then $S=\hat{K}$ is 
a sphere which has the property that $\hat{S}=K$. }

\begin{proof}
If $K=\{x \; \}$, then $\hat{K}$ is the unit sphere $S(x)$.
We know that $\hat{S}$ contains $K$. Assume it contains $K$ and an additional vertex $y$
then $S$ is also the unit sphere of $y$. But then there is an involution which maps $x \to y$ 
and $y \to x$ and leaves $S$ invariant. In other words, $G$ would be reducible and so be not
generic. In general, the fact follows by the recursive definition of genericity. 
If $K=K_{n+1}$, take a vertex $x$ in $K$ and look at the unit sphere. 
It is a sphere of dimension $1$ less and the dual of $\hat{K}$ in $G$ is the
dual of $\hat{K} \setminus \{x\}$ in $S(x)$. By induction, its dual is $K$. 
\end{proof}

The sphere bundle over the simplex set comes handy when filling the interior
structure of the faces of the dual graph $\hat{G}$ of $G \in \Gcal_d$. This is needed 
when we want to get the completion $\overline{G}$ as a graph in $\Gcal_d$ with
dual Betti data. For $G \in \Gcal$ the sphere bundle reveals the dual geometric graph $\overline{G}$ of $G$
as probably Whitney first made clear, of course using  different language. We want to see this
graph theoretically. The spheres in the sphere bundle have dimension $\leq d-1$.
It is the dual graph $\hat{G}$ which is the skeleton containing the $d$-dimensional simplices  of $G$
as vertices. The sphere bundle helps to fill this skeleton with a triangulation
in order to get $\overline{G}$ as a geometric graph in $\Gcal_d$. First of all we need additional
vertices in $\overline{G}$. Take for that the set of spheres $S(x)$. The dual graph $\overline{G}$ has now $v_0+v_d$
vertices. Connect a vertex $x$ with all maximal simplices which contain it. This produces $V_{d-1}(x)$
edges additionally to the already existing ones in $\hat{G}$. We see that these new edges
correspond to vertices in the dual graph of $\hat{S(x)}$. Now we have to complewte the dual
graph of $S(x)$. But that is a lower dimensional problem which shows that all the filling can be done
inductively. In two dimensions, the task was so easy because the dual graph of $S(x)$ is again a circular graph 
which does not need to be completed.  \\

If $G \in \Gcal_d$ is orientable, then $\overline{G}$ is again a graph in $\Gcal_d$ and 
Poincar\'e duality holds. For example, the dual graph of the octahedron is the cube,
a one-dimensional graph with square faces which can be caped in a compatible way.
For a non-orientable graph $G \in \Gcal_2$ but we can not fill the interiors 
in a compatible way to get a dual completion $\overline{G}$ in $\Gcal_2$. To prove the
Poincar\'e duality, one needs to write down the incidence matrices $\hat{d}_k$ of the dual graph
and relate the kernel of the Laplacians $\hat{L}_l = (\hat{D}^2)_l$ with the kernel of $L_{d-k}$
which is not obvious as the matrices have completely different size. The new $\hat{L}_0$ for
example is $(v_0 + v_d) \times (v_0 + v_d)$ matrix. We have not shown yet Poincar\'e duality purely
graph theoretically. 

\resultlemma{
$G$ is Eulerian if and only the completed dual $\overline{G}$ is Eulerian.
}

\begin{proof}
Proof by induction for $d$. For $d=1$, it is clear as the dual graph is isomorphig to the graph itself. 
Look at the unit spheres in the completed graph. For the $v_d$
vertices belonging to the maximal simplices, it is the completion of the dual graph 
of the sphere $S(x)$. By induction and using that $S(x)$ is Eulerian, also the dual
is Eulerian if and only if $S(x)$ is Eulerian. For each of the $v_0$ virtual vertices $y$ added, 
the unit sphere $S(y)$ is isomorphic to $S(x)$ and so Eulerian.
\end{proof} 

{\bf Examples.} \\
{\bf 1)} For $d=1$, we always have $\overline{G}=G$.  \\
{\bf 2)} For $d=2$, the volumes of the graph $\overline{G}$ are $\hat{v}_0 = v_0 + v_2$
$\hat{v}_1 = 2 v_1$ and $\hat{v}_2 = v_1$. The graph $\overline{G}$ 
now has vertex degrees $6$ or $V_0(x)$, which are all even. \\
{\bf 3)} For the 16-cell $G$, the dual graph is the stellated three dimensional cube
which is Eulerian as every unit sphere is either an octahedron or stellated cube. \\
{\bf 4)} We see already when looking at the two dimensional spheres that when 
taking the double dual, the volumes grow exponentially. \\

\resultlemma{
Taking completed double duals is a way to refine a graph in a way that
$\Ecal_d$ is left invariant. 
}

\begin{figure}[h]
\scalebox{0.36}{\includegraphics{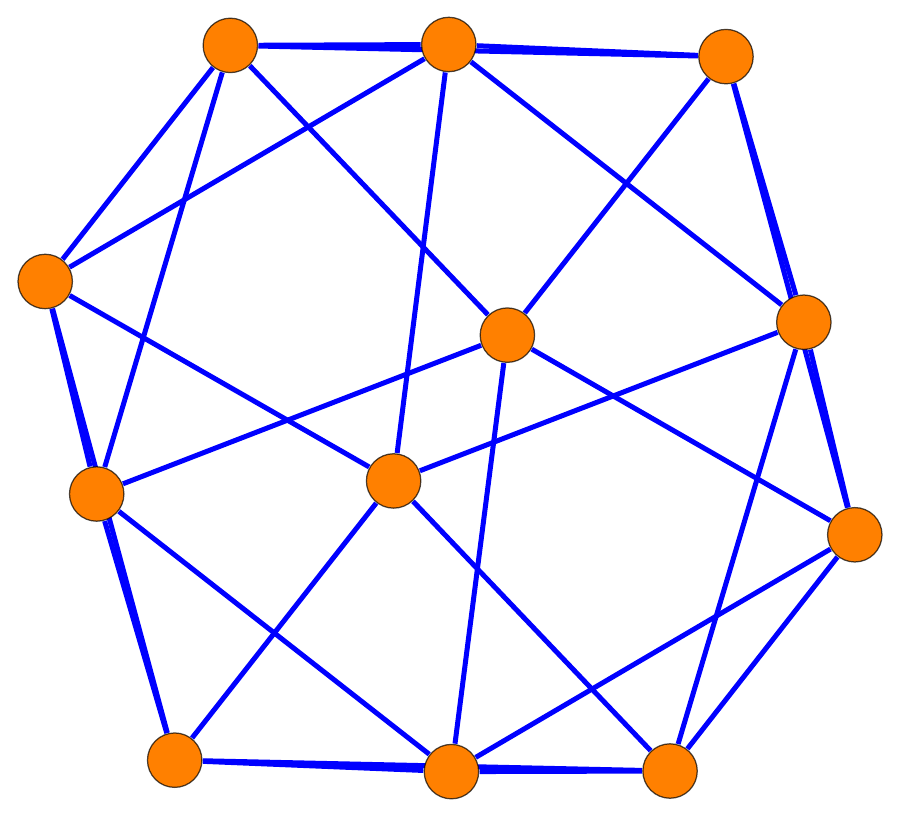}}
\scalebox{0.36}{\includegraphics{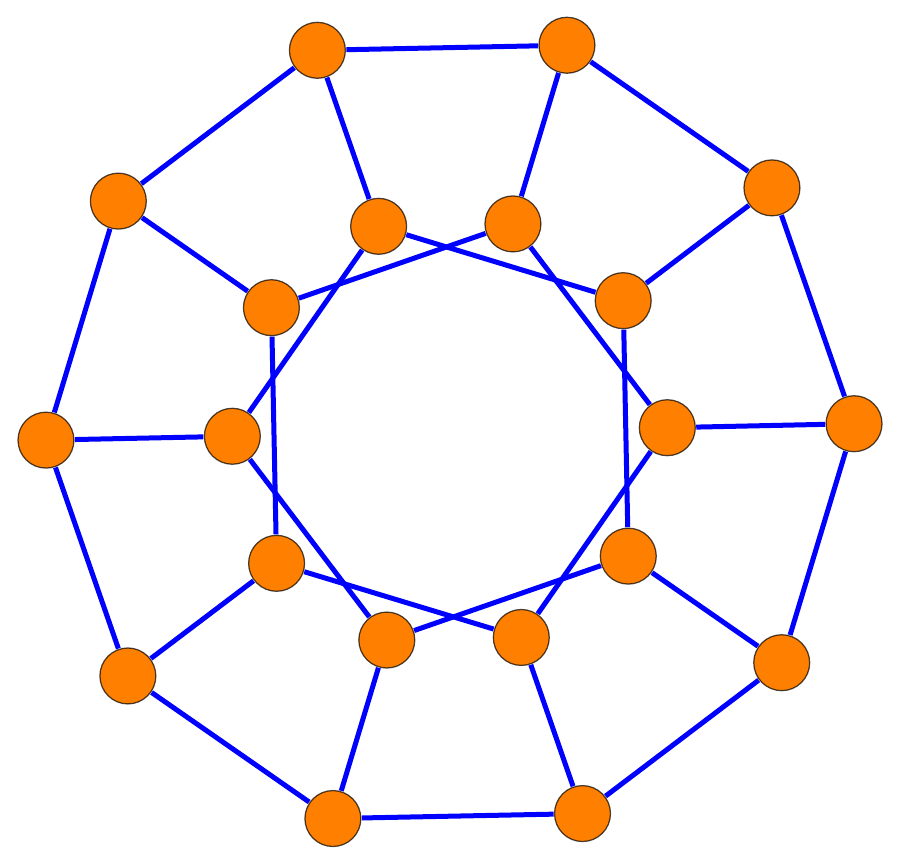}}
\caption{
The icosahedron graph $G$ and its dual, the dodecahedron $\hat{G}$. The dual 
graph is not bipartite which is always the case if $G \notin \Ecal_d$. 
The icosahedron is a prototype graph, where we have no natural geodesic flow
as at each cross road we have to make a choice which opposite road to chose. 
}
\end{figure}

\begin{figure}[h]
\scalebox{0.36}{\includegraphics{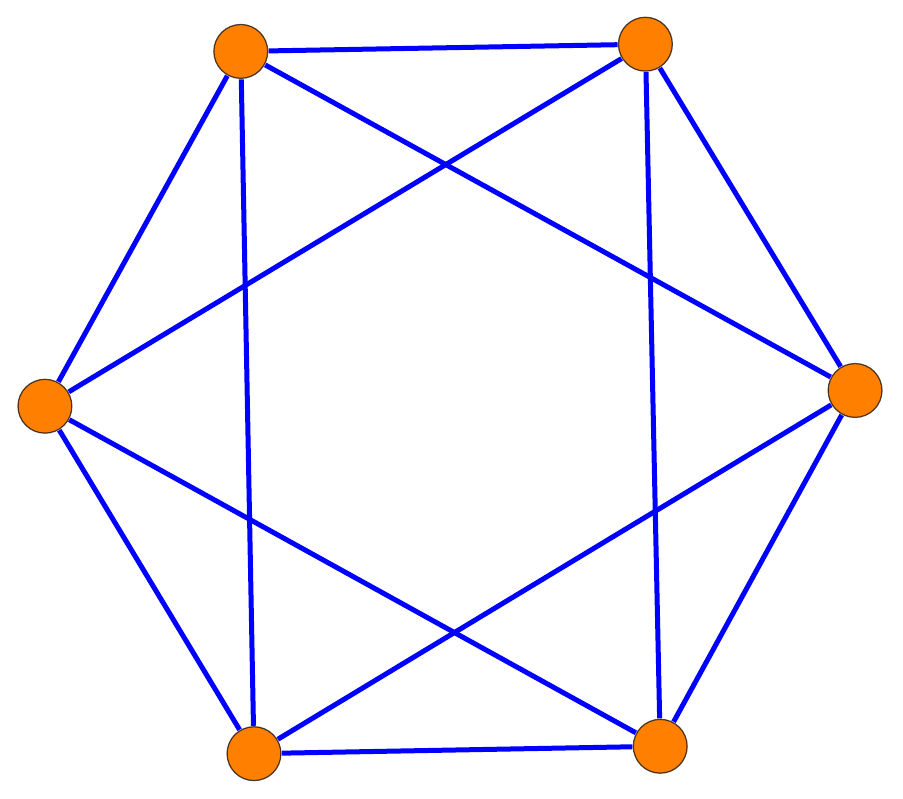}}
\scalebox{0.36}{\includegraphics{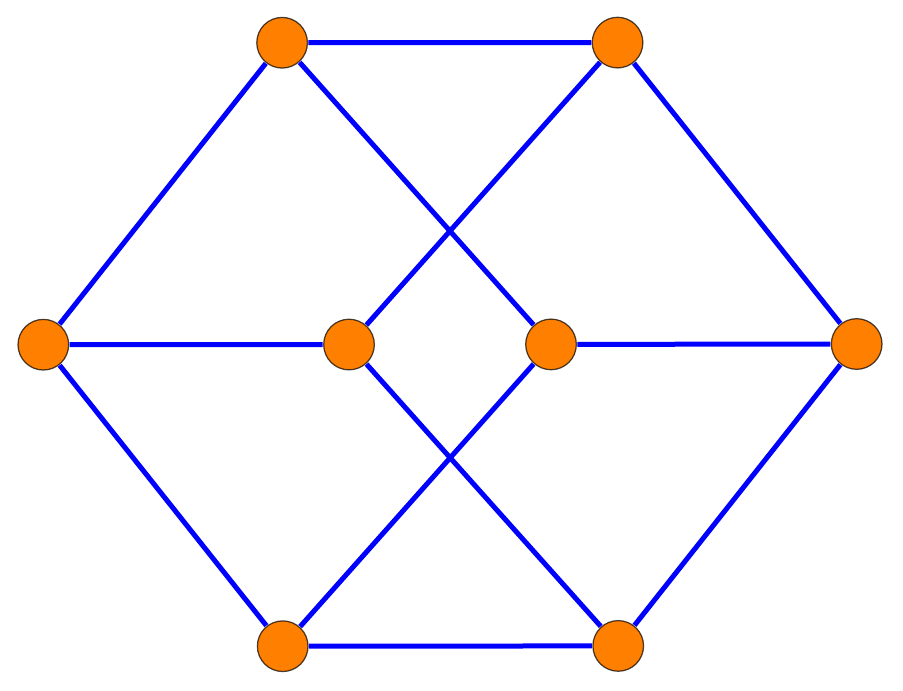}}
\caption{
The octahedron graph $G$ and its dual, the cube $\hat{G}$. The dual 
graph is bipartite which is always the case if $G \in \Ecal_d$. 
The octahedron is a Platonic sphere. All Platonic spheres in $d=4$ and higher
dimensions can be drawn by taking a $2(d+1)$-gon and draw all diagonals except
opposite diagonals. 
}
\end{figure} 

\section{Sphere bundles over the simplex graph} 

\begin{figure}[h]
\scalebox{0.36}{\includegraphics{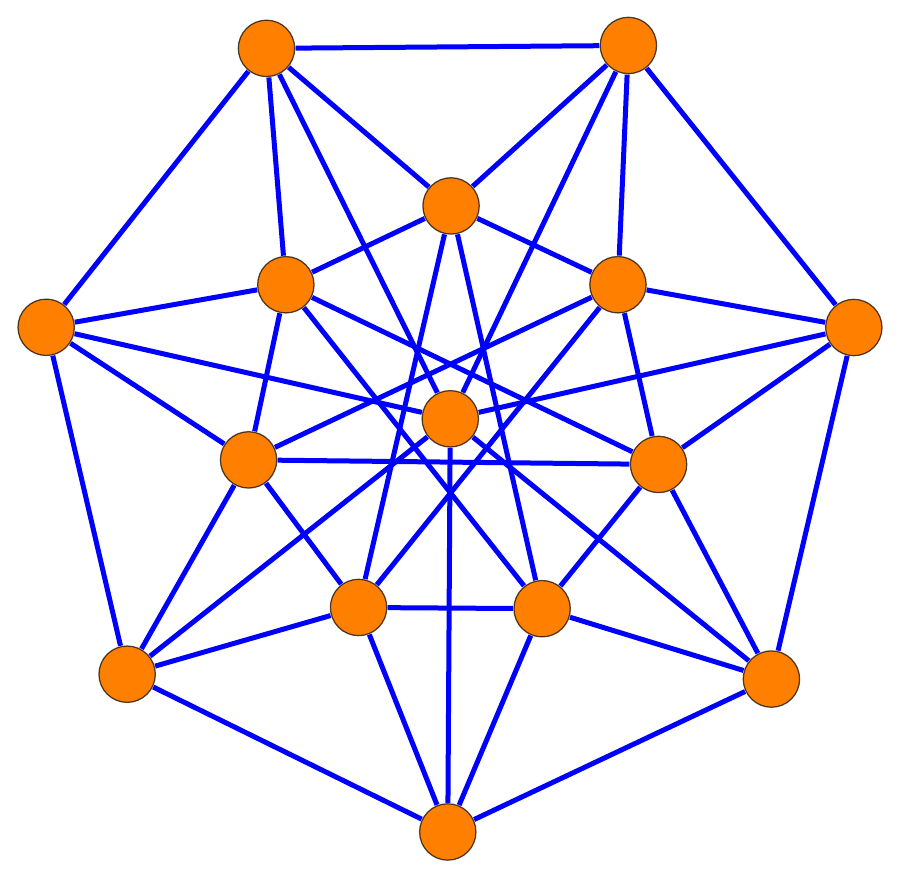}}
\scalebox{0.36}{\includegraphics{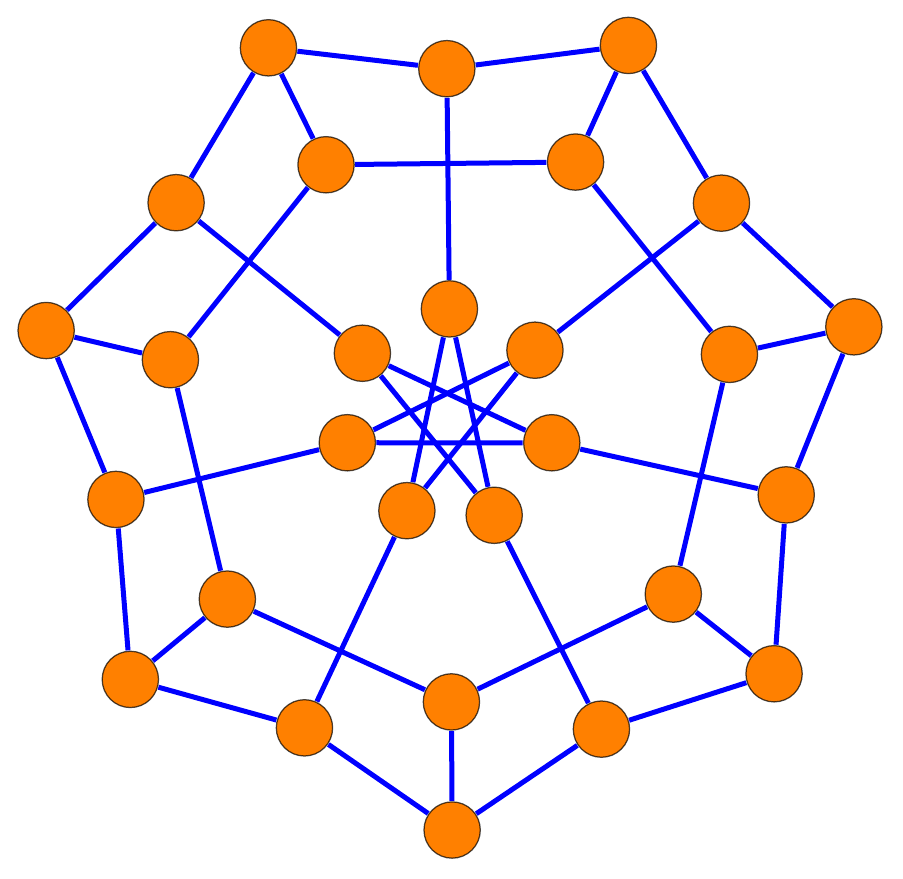}}
\caption{
Jenny's graph $G$, a projective plane and its dual. This graph $G \in \Gcal_2$ 
and has chromatic number $5$. The graph has been found in the Summer 2014 by Jenny Nitishinskaya
during a HCRP project.  It is non orientable and the bipartite structure 
says nothing about the chromatic number. 
There are projective planes with chromatic number $3$ for which the 
dual graph is still not bipartite. 
}
\end{figure}

Simplices play an important role in graph theory. In geometric setups, this has been
realized most prominently by Forman \cite{forman95} who built discrete Morse theory. 
Calculus and cohomology are other reasons as differential forms are just functions
on simplices in such a way that an pre-described ordering of the vertices in each 
simplex produces functions on the simplices. 
An other indication comes from fixed point results like the Lefshetz fixed point 
formula \cite{brouwergraph} which equates the sum of the 
degrees of fixed simplices with the Lefshetz number $L(T)$ of the graph automorphism $T$. 
This result holds in full generality for all finite simple graphs. 
In a geometric setup, if $G \in \Gcal_d$, then the degree of the map $T$ on a simplex $x$ 
can be seen in a more general way as we have now a natural sphere bundle on the set of simplices.
We have just seen that under some genericity condition, the sphere recovers the simplex. \\

Assume now that $G$ is in $\Gcal_d$ and that $S(x)$ denotes the sphere associated to a 
complete subgraph $x$ which we consider as in the Brouwer story to be an atomic quantity. 

\resultlemma{
If a graph automorphism $T$ has a fixed simplex $x$ then 
it induces an automorphism on its dual sphere $\hat{x}$. }

\begin{proof}
As $T(x)=x$, the permuted $x$ has the same sphere.
Let $a_1,\dots ,a_k$ are vertices in $x$ which are cyclically permuted, 
this cyclically permutes the corresponding spheres $S(a_i)$, leaving the 
intersection invariant. It therefore produces a map on the sphere $\hat{x}$. 
As $T$ is an automorphism, it induces an automorphism on the subgraph $\hat{x}$. 
\end{proof}

{\bf Examples.} \\
{\bf 1)} If $G$ is the octahedron and $T$ is a rotation of this sphere with $2$ fixed 
vertices, then these are the only simplices $x,y$ fixed and their degree is one. The degree of
the rotation on $\hat{x}=\hat{y}$ is $0$. 
The Lefshetz number of $T$ on $G$ is $2$ as the map $T$ induces the identity on 
$H^0(G)$ and $H^2(G)$ and has no contribution from $H^1(G)$ as $G$ is simply connected. \\
{\bf 2)} If $G$ is the octahedron and $T$ is the antipodal map, then $T$ is 
orientation reversing and has trace $-1$ on $H^2(G)$. Therefore, its Lefshetz number is zero. 
There is no fixed point. \\

This leads to a reformulation of the Brouwer-Lefshetz story as if $x$ is a fixed simplex, then
the dual graph $\hat{x}$ is fixed too.

\definition{ If $T$ is an automorphism of a geometric graph, 
define the degree of $T$ on $\hat{x}$ as the degree of $x$. }

The Brouwer-Lefshetz fixed point theorem \cite{brouwergraph} can now be formulated in that the sum of the degrees
of $T$ induced on spheres is the Lefshetz number. In some sense this is closer to the 
continuum geometry, but its not clear yet how things are linked. The Lefshetz number of 
a sphere is either $2$ or $0$ depending on the dimension and whether $T$ is orientation preserving
or not. If $d-k-1$ is even and $T$ is orientation preserving on the sphere $\hat{x}$ or $d-k-1$ is odd and $T$ is orientation 
reversing on the sphere $\hat{x}$ then $L(T|\hat{x})-1 = 1$, otherwise, it is $-1$. 
The degree of $x$ is $1$ if $k$ is even and $T$ is orientation preserving on $x$ or if $k$ is odd and 
$T$ is orientation reversing on $x$. We see that the situation is completely parallel. \\

It is surprisingly tricky to define without intersecting spheres what a 
``lower dimensional unit sphere" is in a graph. In the case of one dimensional spheres,
we can not just say it to be a subgraph isomorphic to 
$\Scal_{1}$ for which the diameter in $G$ is $2$ as any closed curve on the unit 
sphere $S(x)$ of a vertex would qualify as a $1$-dimensional unit sphere. 
Duality provided us with a convenient way:  $k$-dimensional sphere in a $k$-dimensional 
geometric graph is the dual graph of a $K_{d-k}$ subgraph. 
For example, for $d=3$, a $k=1$ dimensional subgraph is the dual 
of a $K_2$, an edge. The {\bf sphere bundle} of the set of simplices is therefore a natural 
construct and produces all the small spheres we ever need in a geometric setup. 

\section{Geodesics and Hopf-Rinov}

There will never be a completely satisfactory classical geodesic flow for finite simple 
graphs simply because the unit sphere at a point $x$ is in general
smaller than the number of vertices different from $x$. Even if a unique geodesic 
flow can be defined in such a way that two points in distance $2$ have a unique geodesic connection,
still many graphs have the property that for some point $x$, the union of all geodesics starting 
at $x$ will not cover the entire graph. The true geodesic flow is {\bf quantum} in nature, as it is in the 
real world: we have to look at the wave equation on a graph. It is then possible for two points $x,y$ 
to start with a wave located on the vertex $u(0)=1_{\{x\}}$ at time $t=0$ and find a velocity vector $u'(0)$ 
of length $1$ (a function on vertices and edges) and a time $T$ such that $u(T) = 1_{\{y\}}$. 
The reason why this is possible is simply linear algebra: because 
the wave equation $u''=-Lu$ with scalar Laplacian of the
graph is part of the wave equation $u''=-Lu$ on the simplex space, where $L=D^2$ with 
a matrix $D=d+d^*$ so that we can write $u(t)=\cos(Dt) u(0) + \sin(Dt) D^{-1} u'(0)$ which is the
d'Alemberg solution. We now only have to solve for $t$ and $u'(0)$ in the ortho-complement of the kernel
of $D$ to solve the equation for given $u(t)$ and $u(0)$. The smallest time $T>0$ which leads to a solution gives
us a notion of ``how fast particles travel" in the graph. 
Mathematically, there is no problem in the simply connected connected case as we only have to 
look for velocities for which the total velocity $\sum u'(0,x)$ is zero, to make sense of$D^{-1} u'(u)$. 
As $D$ is a symmetric matrix, there is no need even to involve the theory of pseudo inverses. In any case,
this solution of the wave equation is quantum in nature because we can write also 
$\psi(t) = e^{i t D} \psi(0)$ with complex wave function vectors 
$\psi(t) = u(t) + i D^{-1} u'(t)$. So, waves on a graph are described by the Schr\"odinger equation
for the Dirac operator.  \\

For a physicist, this story really becomes exciting when looking at the
propagation of waves on higher forms. For example, if $u(0)$ is an initial wave located on an edge,
then we look at the propagation of $1$-forms which is electromagnetism. For $2-$forms given by functions
on triangles, the symmetry group of the triangle is already non-Abelian. How fast do things travel there
answers the question ``how heavy" the particles under consideration are. In any case,
the mathematics of propagating waves (= quantum particles) in this finite universe is simple 
linear algebra even so the matrices can become large. For the icosahedron already, the matrix $D$
is a $62 \times 62$ matrix as $v_0+v_1+v_2=12+30+20=62$. For a three dimensional graph modeling 
more realistic physics, the matrices are much larger. \\

If we go back to draw circles and look for a classical notion of ``line" we need to say what a geodesic is
in the graph. We can look for a weaker form of Hopf-Rinov and not ask to be able to connect two
arbitrary points with a line - which as we have just seen is impossible - but to extend a line 
segment of length one indefinitely. In other words, we only want to find a dynamical system on the 
unit sphere bundle which has global existence of solutions. Here is where the relation with
chromatology comes in. If we look at a polyhedron $G \in \Scal_2$ which has a vertex with an 
odd number of edges attached and a light ray comes in, then how do we tell, where it goes out? 
We could chose randomly but we can never do that in a reversible manner except bouncing
back some light ray from one direction which is highly unphysical and contradicts anything
we know about light.  Mathematically, the reason why we can not find a local deterministic
reversible propagation rule on a unit ball of an odd degree vertex is that on $C_n$ with odd $n$ 
there is no fixed point free involution. So, we need each vertex to have even degree. That means
the polyhedron has to be Eulerian or equivalently to be $3$ colorable. In higher dimensions,
where unit spheres have dimension larger than $2$, the unit spheres can become already complicated
and we need a bit more than the graph just to be Eulerian. This is why we have looked at 
projective graphs. This is a natural as for large enough unit spheres, the notion of projective
means that the unit sphere is a double cover of a projective space as in Euclidean geometry. \\

\begin{figure}[h]
\parbox{15cm}{
\parbox{6.5cm}{\scalebox{0.15}{\includegraphics{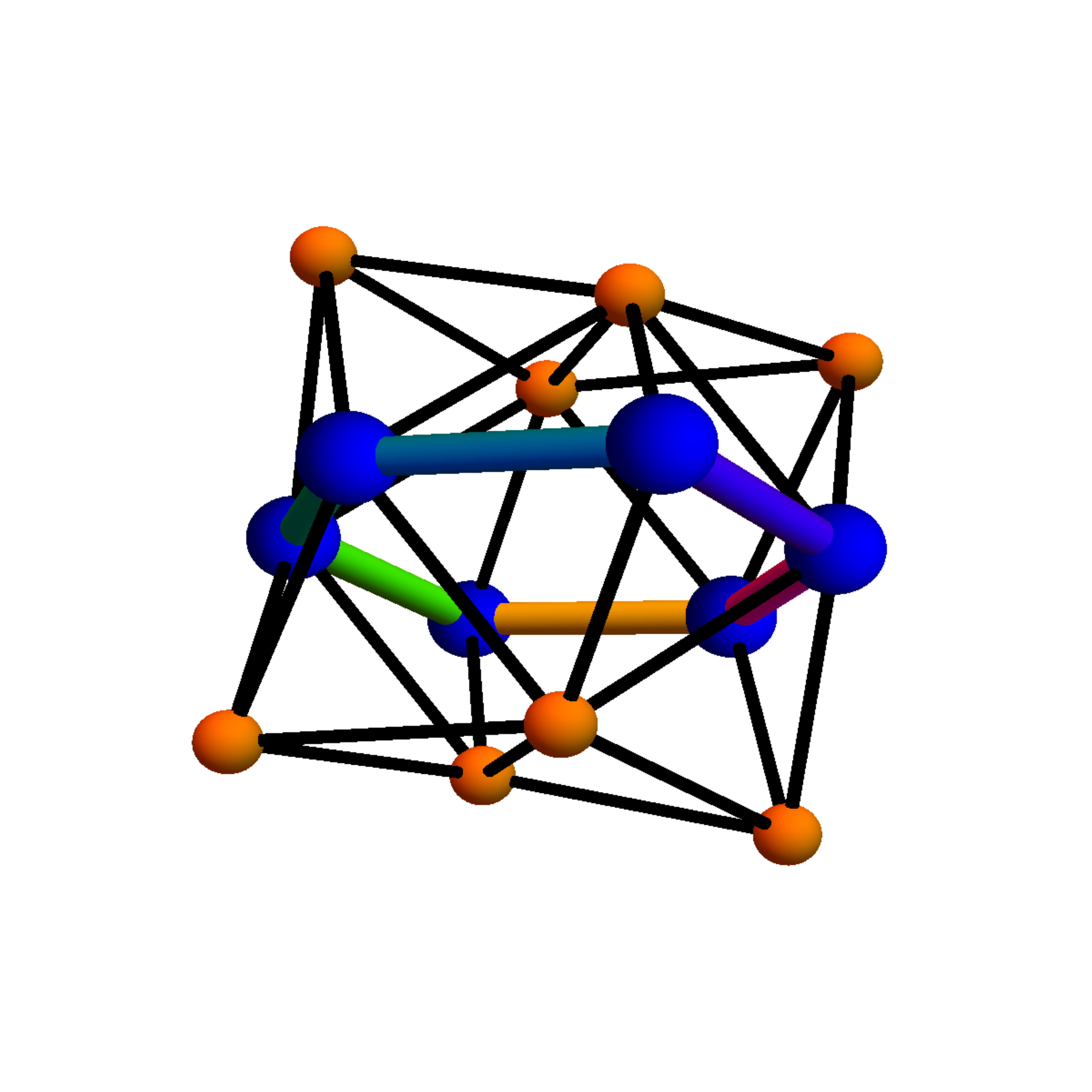}}}
\parbox{6.5cm}{\scalebox{0.15}{\includegraphics{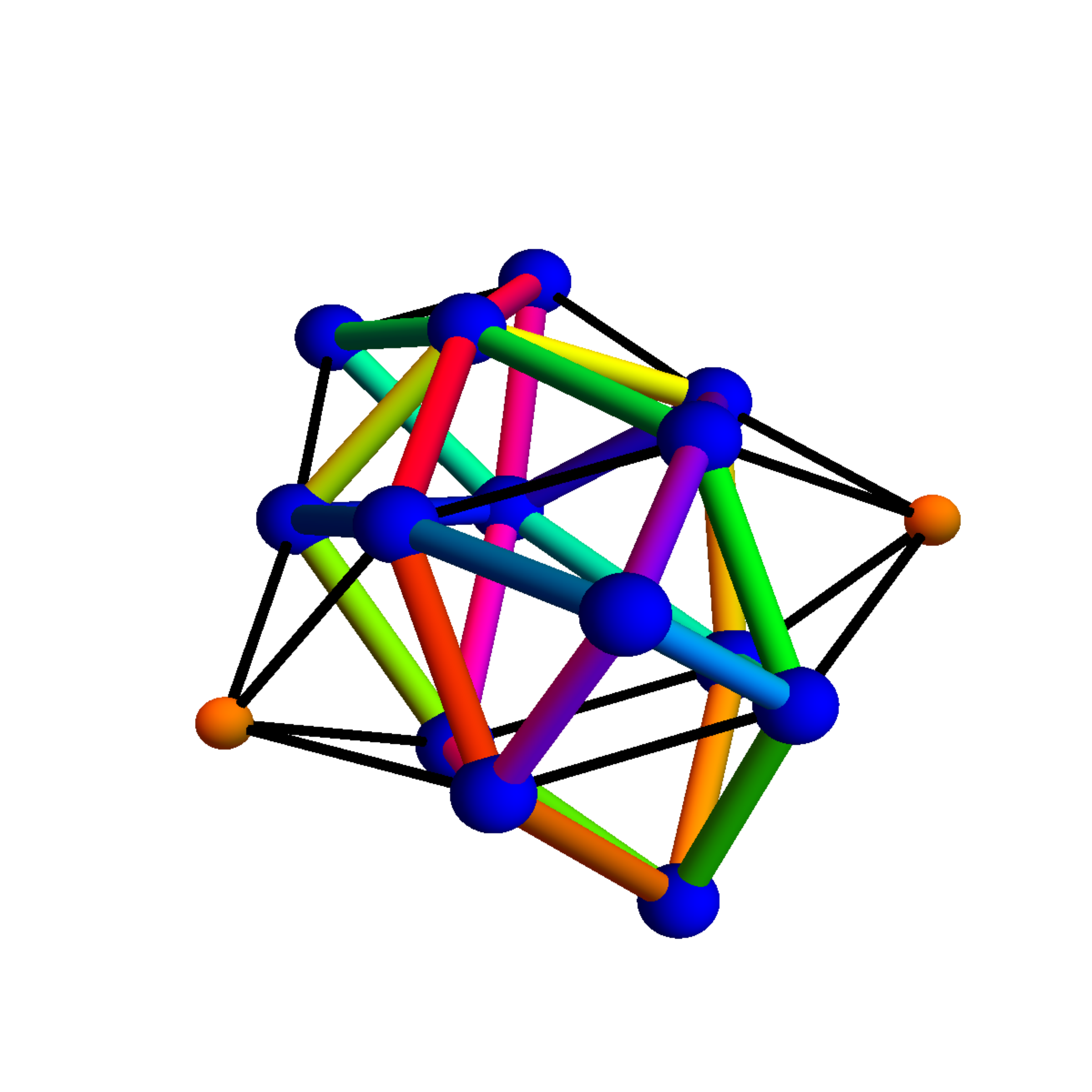}}}
}
\caption{
We first see a closed geodesic path in an Eulerian graph. All 
geodesics are simple closed curves, a classical problem in differential geometry
asks for such manifolds. We see in this case already that the exponential map
is not surjective at some vertices. In the second figure, the graph has 
been refined a bit It had been obtained from 
a capped cube by a double edge refinement so that it is again Eulerian.  }
\end{figure}

The problem of computing geodesic paths in graphs was our very first 
entry point to the field of graph theory while investigating 
classical geodesic evolutions in a HCRP project with  
Michael Teodorescu in the fall of 2008 and spring 2009, especially in the 
context of the open Jacobi conjecture about caustics in ellipsoids.  \\

In any case, the notion of projective spheres is crucial for getting graphs which have geodesic flow 
without having to refer to the quantum world, where waves are computed on graphs. 
For geometric purposes, it is desirable to have a classical Hopf-Rinov type statement. And there
is almost nothing to show: 

\resultlemma{
If a geometric graph $G$ is projective, then there is a unique geodesic flow
on the unit sphere bundle of the graph. 
}

\begin{proof}
The involution $T$ in the unit sphere $S(x)$ describes the propagation of the path in each unit ball:
an incoming ray $(y,x)$ is propagated to a well defined outgoing ray $(x,T(y))$. 
\end{proof}

This is natural since if there exists a fixed point free involution on a sphere, it is unique
among this kind. Comparing with differential geometry, it is a Levi-Civita statement telling
that under some conditions, there is a unique connection, or notion of parallel transport:

\resultlemma{
For a projective sphere $G$ and all simplices $x$, the sphere $S(x)$ is projective. The involution
$T$ on $G$ is unique. 
}
\begin{proof}
Use induction with respect to $d$. 
For $d=1$, it is clear as the automorphism group of a circular graph is the 
dihedral group containing only one involution which is fixed point free and that is the antipodal
map. Assume $T(x)=y$, since $T$ is by induction uniquely determined on $S(x)$, we can extend 
the map uniquely to the disc $B(x)$ from there to a larger neighborhood etc until it fixes the map
on all the sphere. 
\end{proof}

\section{Billiards and caustics}

There is therefore exactly one geodesic flow on a projective geometric graph. The 
projective condition plays the analogue role of zero-torsion for the connection as 
the later also assures uniqueness of the connection. 
And so, we can get geometric graphs with the Hopf-Rinov property if the unit 
spheres are double covers of Eulerian projective spaces. The graph $G$ itself does not
have to have this property. \\

\begin{figure}[h]
\parbox{15cm}{
\parbox{6.5cm}{\scalebox{0.15}{\includegraphics{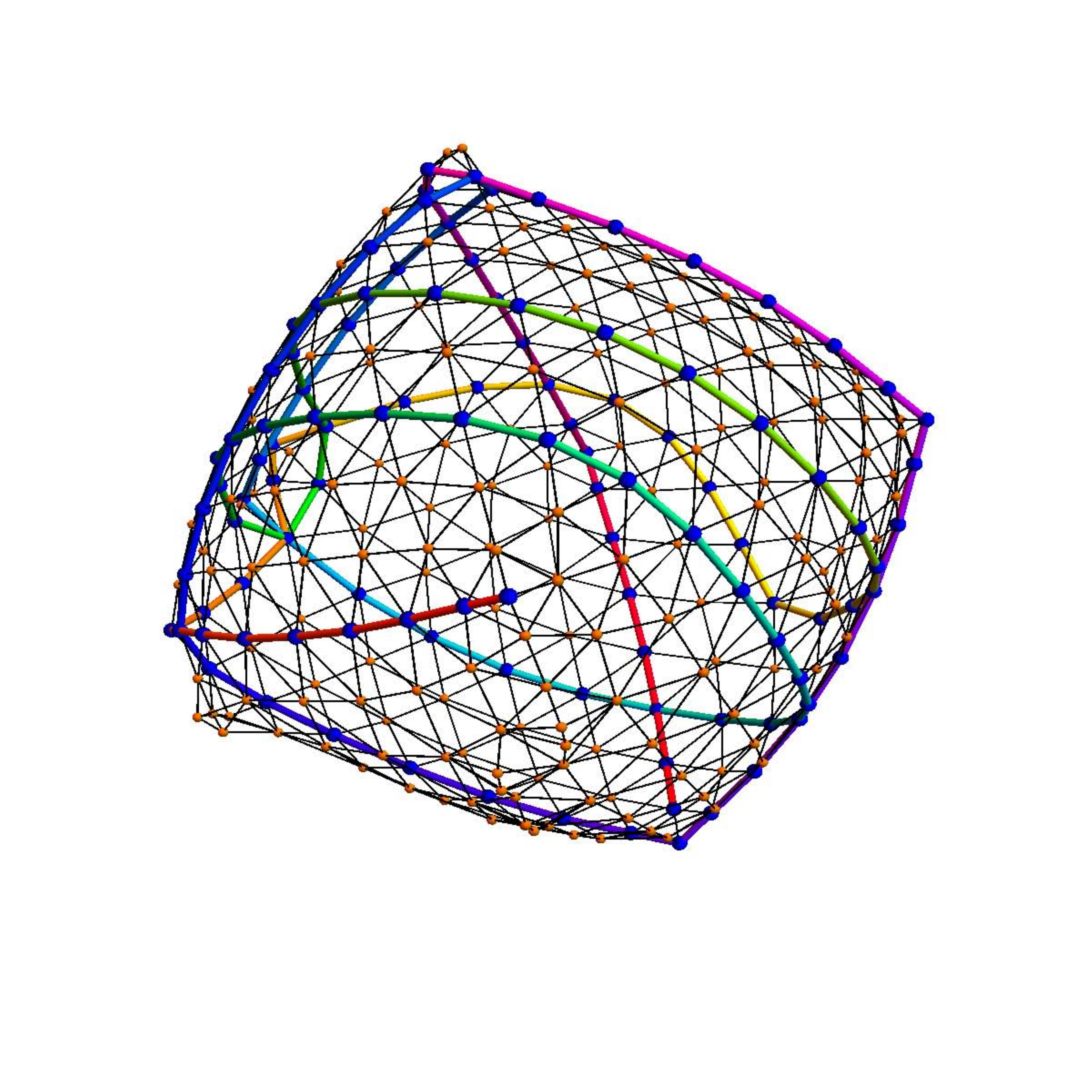}}}
\parbox{6.5cm}{\scalebox{0.15}{\includegraphics{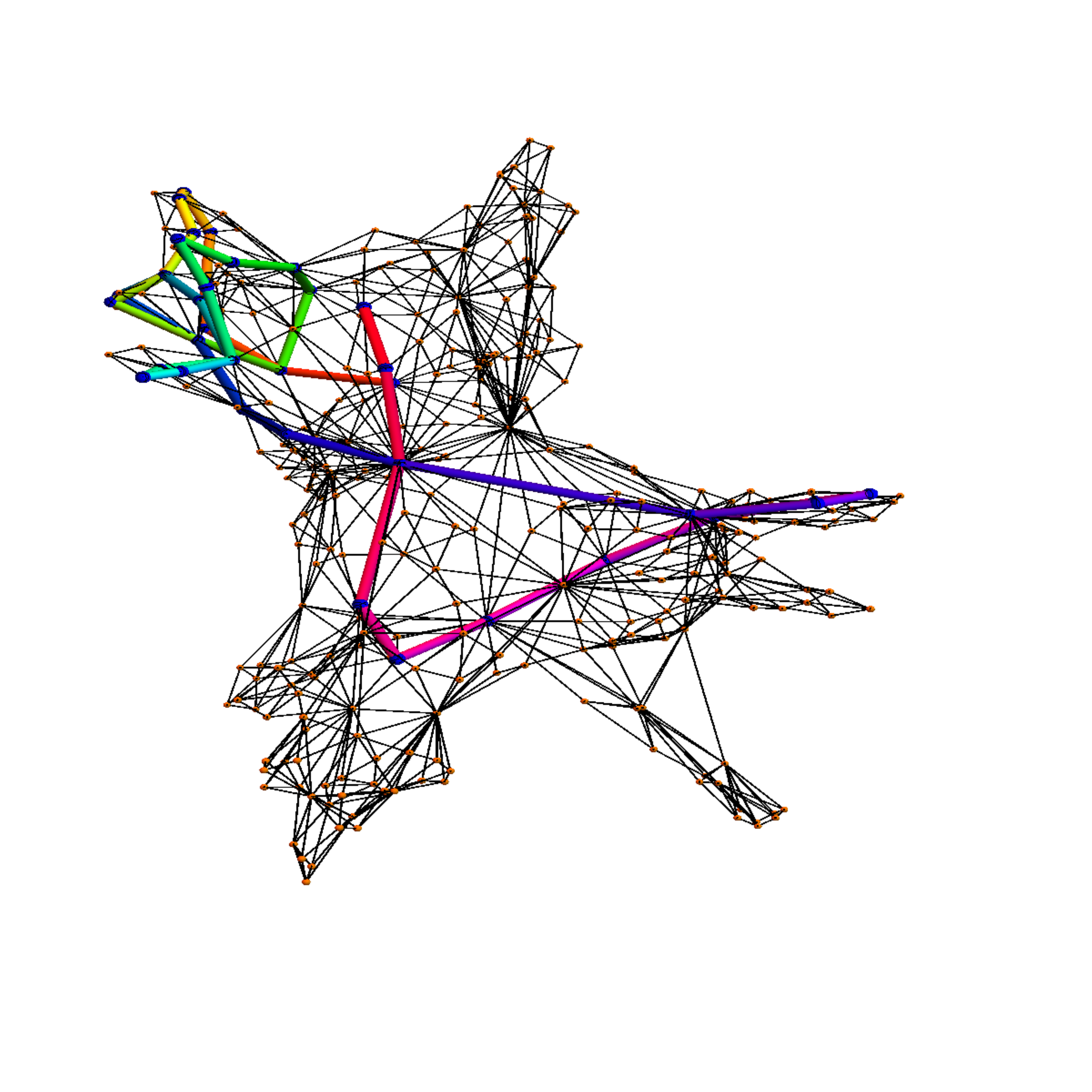}}}
}
\caption{
The first figure shows an Eulerian sphere obtained by subsequent refinements of an octahedron, additionally 
``perturbed" by a couple of edge subdivisions. We see a geodesic of length $140$. 
The second figure shows an Eulerian sphere in $\Ecal_2$ obtained by random subdivisions of a cube graph. 
It looks pretty wild but has dimension 2, Euler characteristic 2 and every vertex
has a circular sphere $S(x)$. In this graph, a geodesic of length 40 is shown. 
}
\end{figure}

\begin{figure}[h]
\parbox{15cm}{
\parbox{6.5cm}{\scalebox{0.15}{\includegraphics{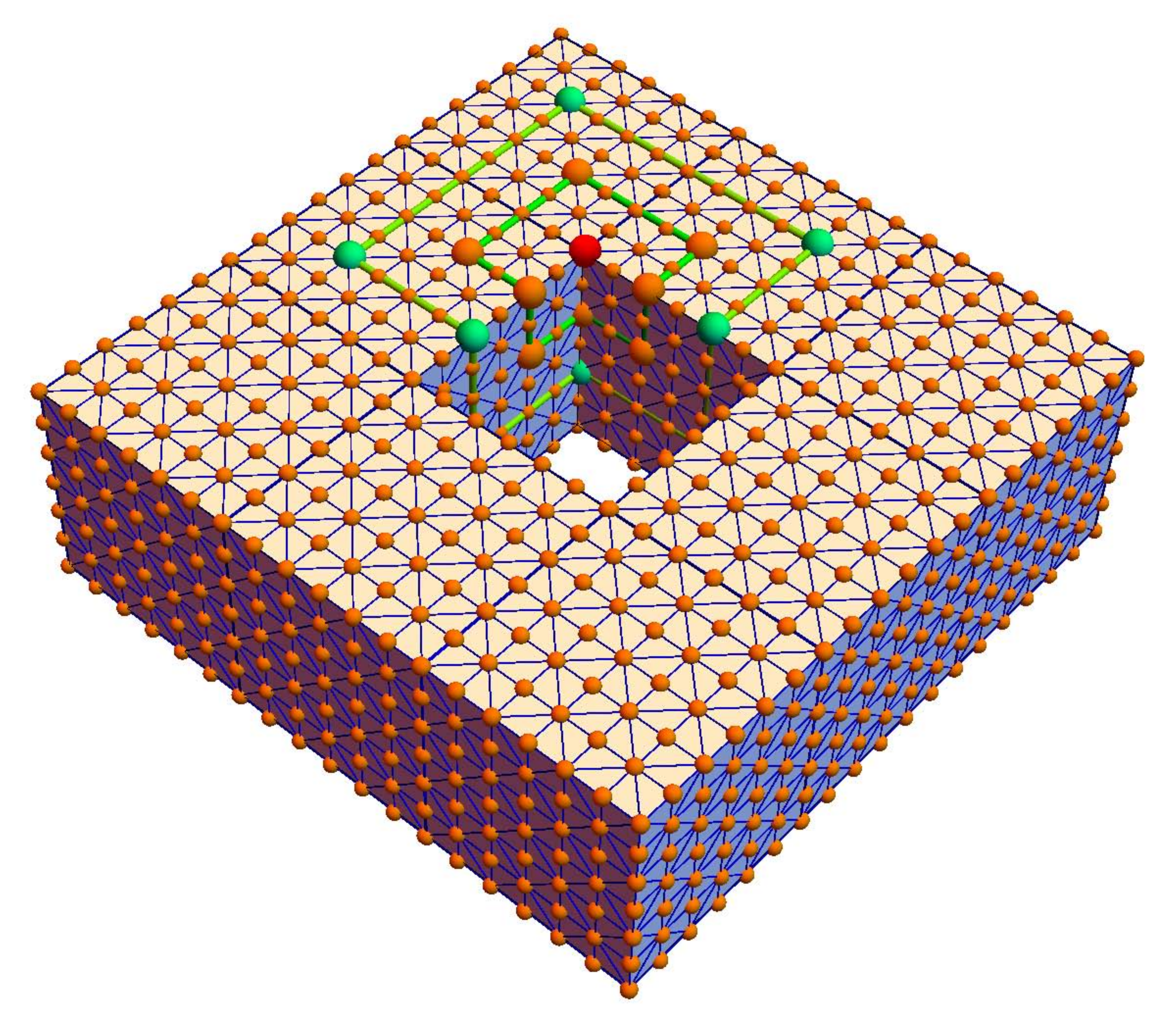}}}
\parbox{6.5cm}{\scalebox{0.15}{\includegraphics{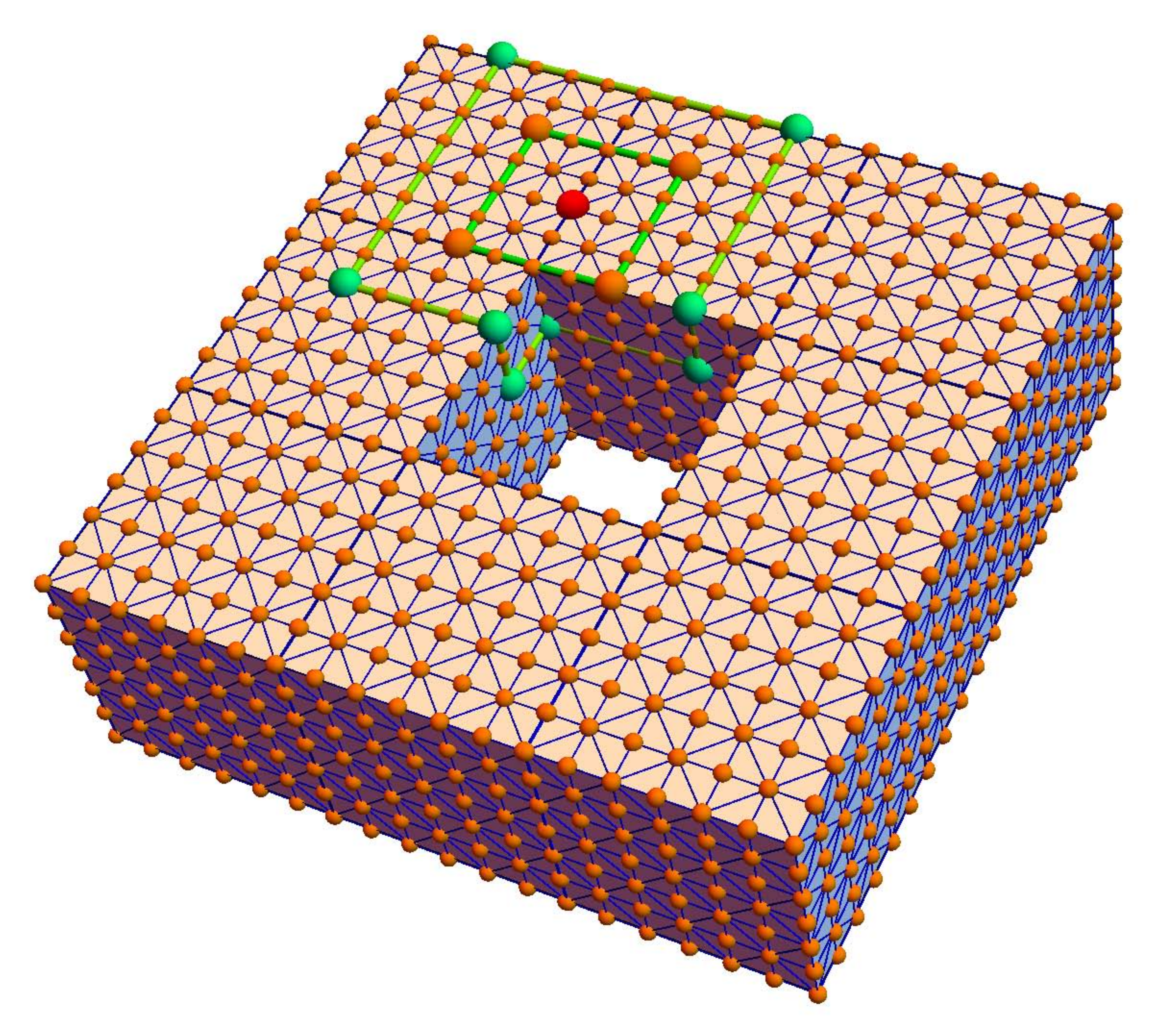}}}
}
\caption{
Wave fronts on an Eulerian graph in $\Gcal_2$. These pictures were
made while investigating the question for which graphs the curvature
$K_2(x) = 2 S_1(x)-S_2(x)$ satisfies Gauss-Bonnet  \cite{elemente11}.
In the new terminology we have there at Eulerian graphs in $\Scal_2$
on the which the geodesic flow is billiards. }
\end{figure}

We have investigated in \cite{elemente11} a very special case of the question: 

\question{
Assume $G \in \Gcal_2$ has the property that $S_2(x) \in \Scal_1$ for all $x$. 
Under which conditions does the second order curvature $K_2(x)=2 S_1(x)-S_2(x)$ 
lead to Gauss-Bonnet $\sum K_2(x) = 30 \chi(G)$? 
}

Here is something we know about ``curvature 60":

\resultlemma{
If $G \in \Gcal_2$ has only degree 5 and 6 vertices and they all have distance
at least $2$ from each other, then $\sum K_2(x) = 30 \chi(G) = 60$.}
\begin{proof}
By the classical Gauss-Bonnet known since the 19th century for polyhedra,
there are exactly 12 degree 5 vertices. Since all neighbors have degree 6, 
the second order curvature is 5 at each of these vertices. For vertices
in distance 2 or larger to the degree 5 vertices, the curvature is 0 as the
2-disk is then completely flat. Also in the immediate neighborhood of a degree $5$
vertex, the curvature is zero, if the vertex has degree 6. 
\end{proof}

We do no know under which conditions such a second order Gauss-Bonnet result holds and therefore
have stuck to first order curvatures, also in higher dimensions. The constant $30$ for the second
order curvature formula is obtained for the icosahedron, where $K_2(x) = 2 \cdot 5 - 5 = 5$ for
all 12 vertices. One problem is to define the second sphere as there are already
points in distance $2$ which can not be reached by geodesics starting at $x$ if the degree of 
$x$ is $4$ and positive curvature implies for Eulerian graphs the presence of degree 4 vertices 
and so large curvature. It might therefore be that we need non-positive curvature for the second
curvature to be a good but this just seems to lead to Buckyball examples like treated in the 
previous lemma. In general, we have had little luck with getting total curvature 60 if there are
negative curvature vertices. In this context there are interesting constraints:
\cite{IKRSS} for example showed recently that (1/6,-1/6) curvature pairs can not be realized alone on
a torus, even so Gauss-Bonnet would allow for that (thanks to  Ivan Izmestiev for sending us
this).  \\

{\bf Examples.} \\
{\bf 1)} For the stellated cube which is a Catalan solid, not all second spheres are circles. The second
order curvatures are $6$ at eight vertices and zero else. The total curvature is 48. \\
{\bf 2)} From the 13 Catalan solids, the duals of the Archimedean solids, there are 4
which are geometric and 3 of them are Eulerian. For none of them, the second spheres are spheres. 
The history of Archimedean solids is a bit murky since Archimedan's own account is lost \cite{coxeter}
but where Archimedes has referred to the cuboctahedron as already been studied by Plato.
The duality might have been first understood in the 14th century by Maurolycus. \\
{\bf 3)} For all ``Buckminsterfullerene" type graphs with degree 5 and 6 vertices
the total second order curvature is 60: the 
proof is that the isosahedral symmetry forces all the 12 degree 5 vertices to be separated apart 
and the second order curvature is 5 in each of these cases and zero else. Examples are ``golf balls". 
Degree 4 vertices do not work as the second spheres $S_r(x)$ are no spheres. In some sense, degree 4
vertices have the effect that the graph is not ``smooth enough" for second order curvatures.

\begin{figure}[h]
\parbox{15cm}{
\parbox{6.5cm}{\scalebox{0.15}{\includegraphics{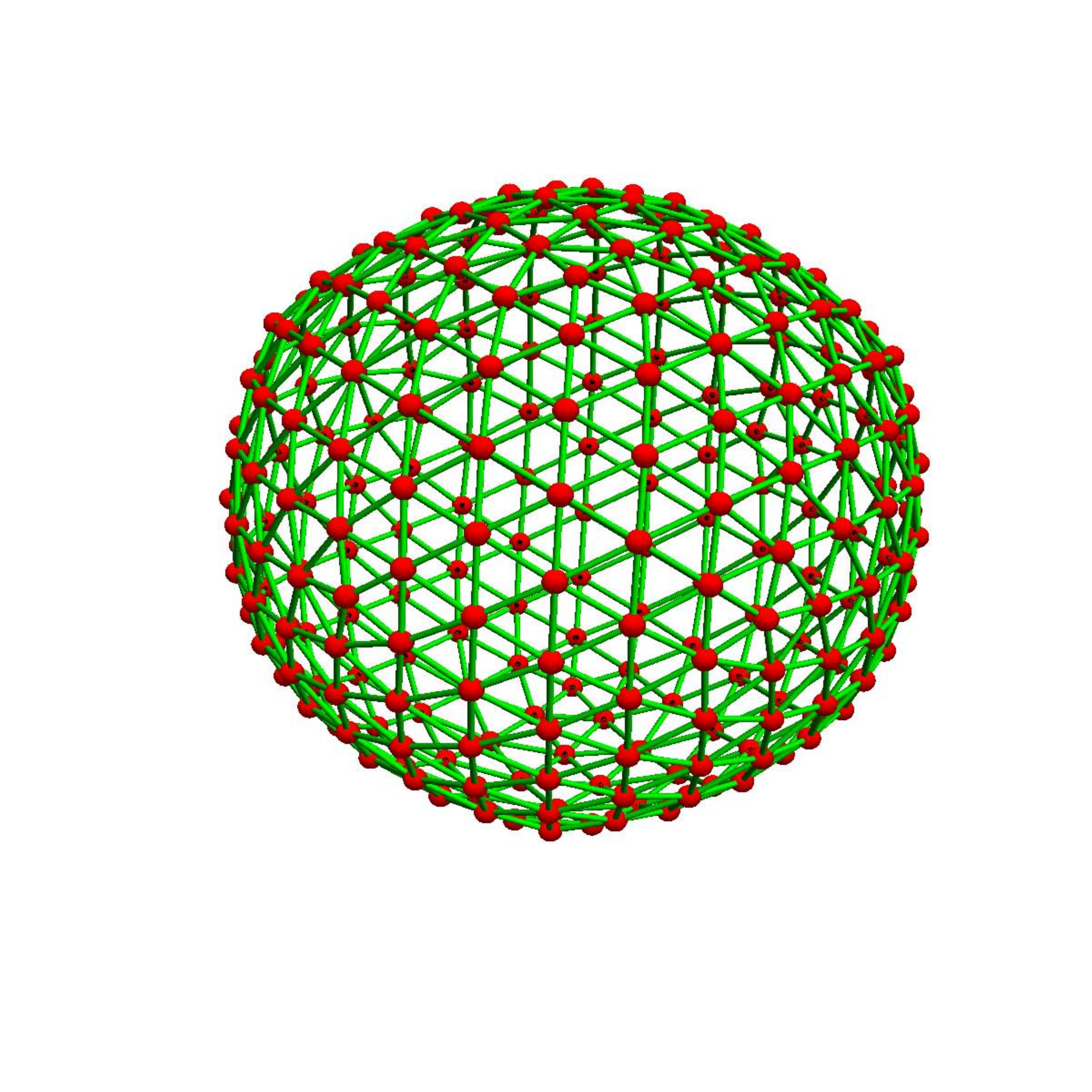}}}
\parbox{6.5cm}{\scalebox{0.15}{\includegraphics{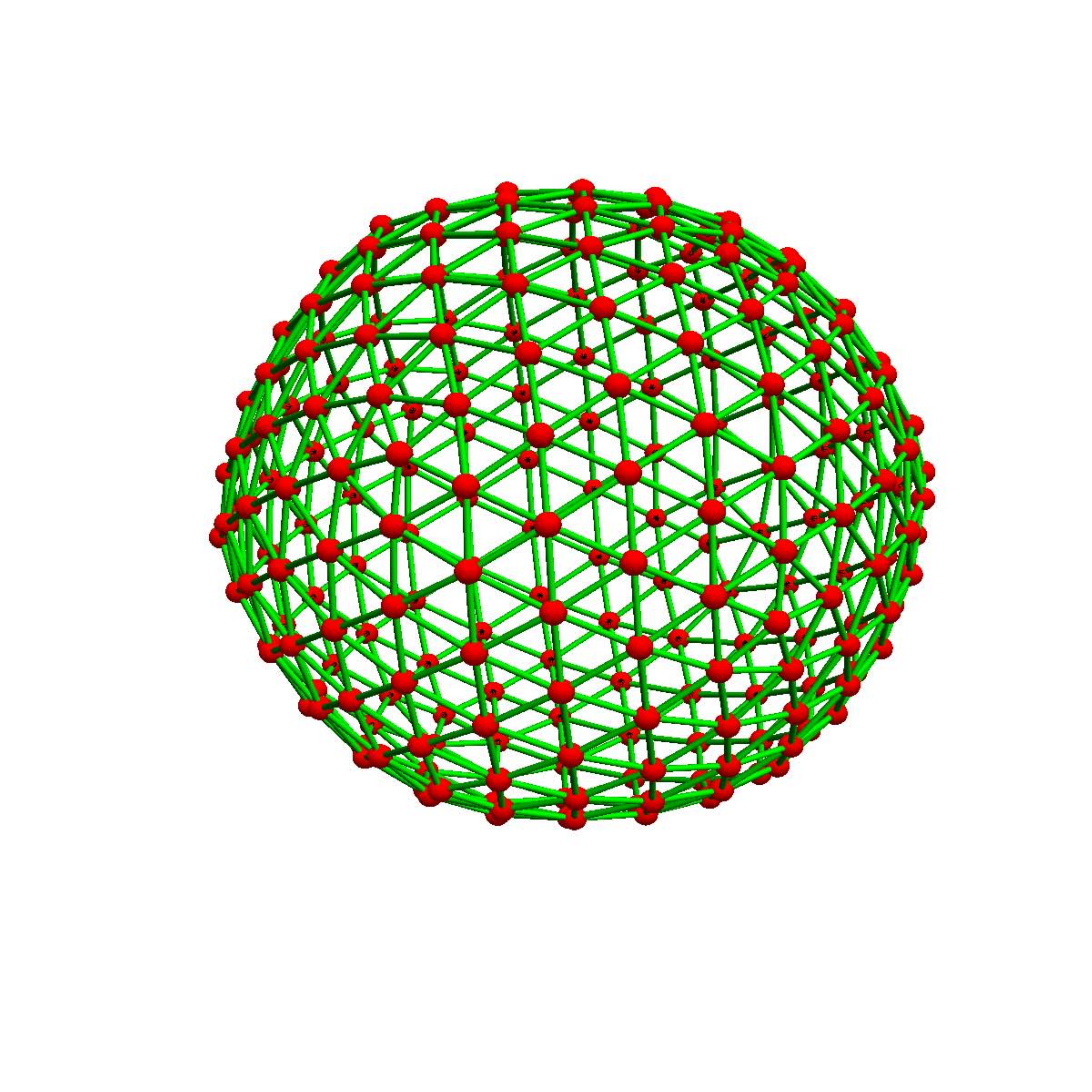}}}
}
\caption{
Buckyballs for which the curvature 60 theorem holds:
the first one has volume data $(v_0,v_1,v_2)=(272,810,540)$
The second one has the volume data $(272,750,500)$. In all cases, the $K_2(x)$
is $5$ at the centers of the 12 degree 5 vertices. 
}
\end{figure}

\begin{figure}[h]
\parbox{15cm}{
\parbox{6.5cm}{\scalebox{0.15}{\includegraphics{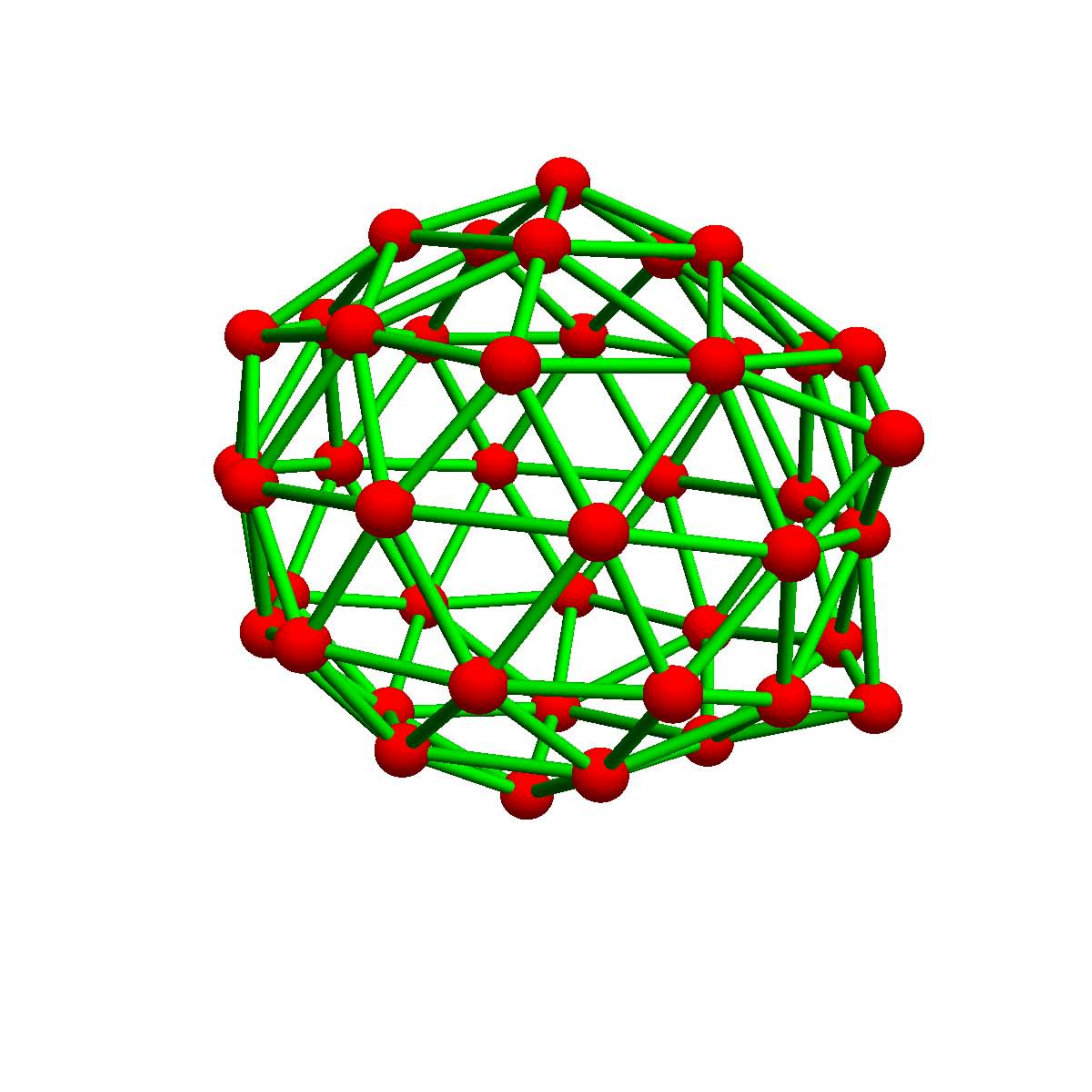}}}
\parbox{6.5cm}{\scalebox{0.15}{\includegraphics{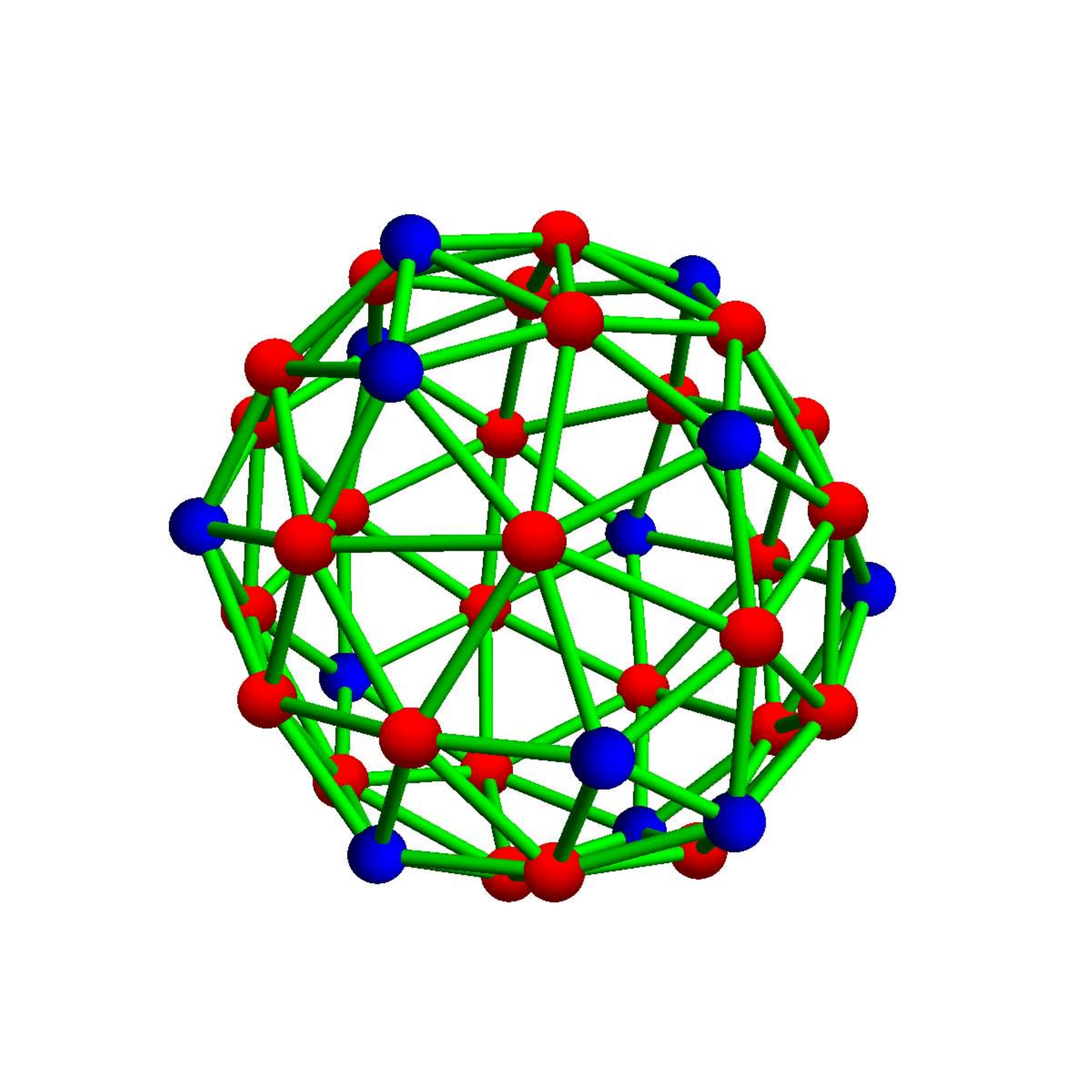}}}
}
\caption{
The first example of a graph obtained by 2 edge subdivisions for which the theorem 60
theorem does not hold. The second order curvatures are 
$(-2, -1, 0, 1, 2, 4, 6)$ which totals to 46, not 60. The volume data of the
polyhedron are $(44,126,84)$. 
The second example is a graph for which the degree is 5 or 6. But there are now 
negative curvature parts. It has been obtained by edge collapse and now there are 
degree 5 leading to a Birkhoff diamond patches and negative curvature points. }
\end{figure}

\definition{
If a projective geometric graph $G \in \Gcal_d$ admits an involution $T$ such 
that a graph $S \in \Gcal_{d-1}$ is fixed, then the geodesic flow is called
{\bf billiards} on $G/T$ with boundary $S$.}

This is a standard motivation for billiards \cite{Moservariations}: classical 
billiards can be seen as a limiting case of a geodesic flow. In the discrete
we do not have to take a limit. The geodesic flow on the double cover of the 
billiard table with boundary is equivalent to the billiard situation, where
the ball changes direction at each boundary point. \\

One can also look at notions of caustic. There are caustics known in billiards
as well as in geodesic flows are sometimes related as \cite{elemente98} illustrates.

\definition{
Let $G \in \Gcal_d$ be a projective geometric graph and $x$ a vertex. The 
{\bf primary caustic} of $x$ is the set of points $y$ for which there
are at least two different geodesics starting at $x$ anding at $y$ and
such that no vertex $z$ on any of the geodesics connecting $x$ and $y$
is already in a caustic.}

A point in the caustic is also called {\bf conjugate point} but it is not the same.
Classically, conjugate points are points for which a nonzero Jacobi field exists along
the geodesics which has a root at the end points. It is a common theme in differential
geometry to give curvature conditions estimating the radius of injectivity in terms 
of curvature (see e.g. \cite{BergerPanorama}). \\

The caustics are difficult to understand in differential geometry is illustrated by the unsolved
{\bf Jacobi's last theorem} asking whether on a general ellipsoid, all primary caustics
have 4 cusps. This is rather embarrassing since the geodesic flow on the ellipsoid is integrable. 
It is the prototype of a system which is ``solvable". Still we have no clue about the caustics. 
One can look at secondary and ternary caustics etc. In general, even in integrable situations, the
union of all caustics is expected to become dense except in very special cases like the round sphere
where the caustic is always a point. It was this story of caustics which brought us originally 
to the journey to consider at discrete versions of the problem.  \\

\begin{figure}[h]
\scalebox{0.30}{\includegraphics{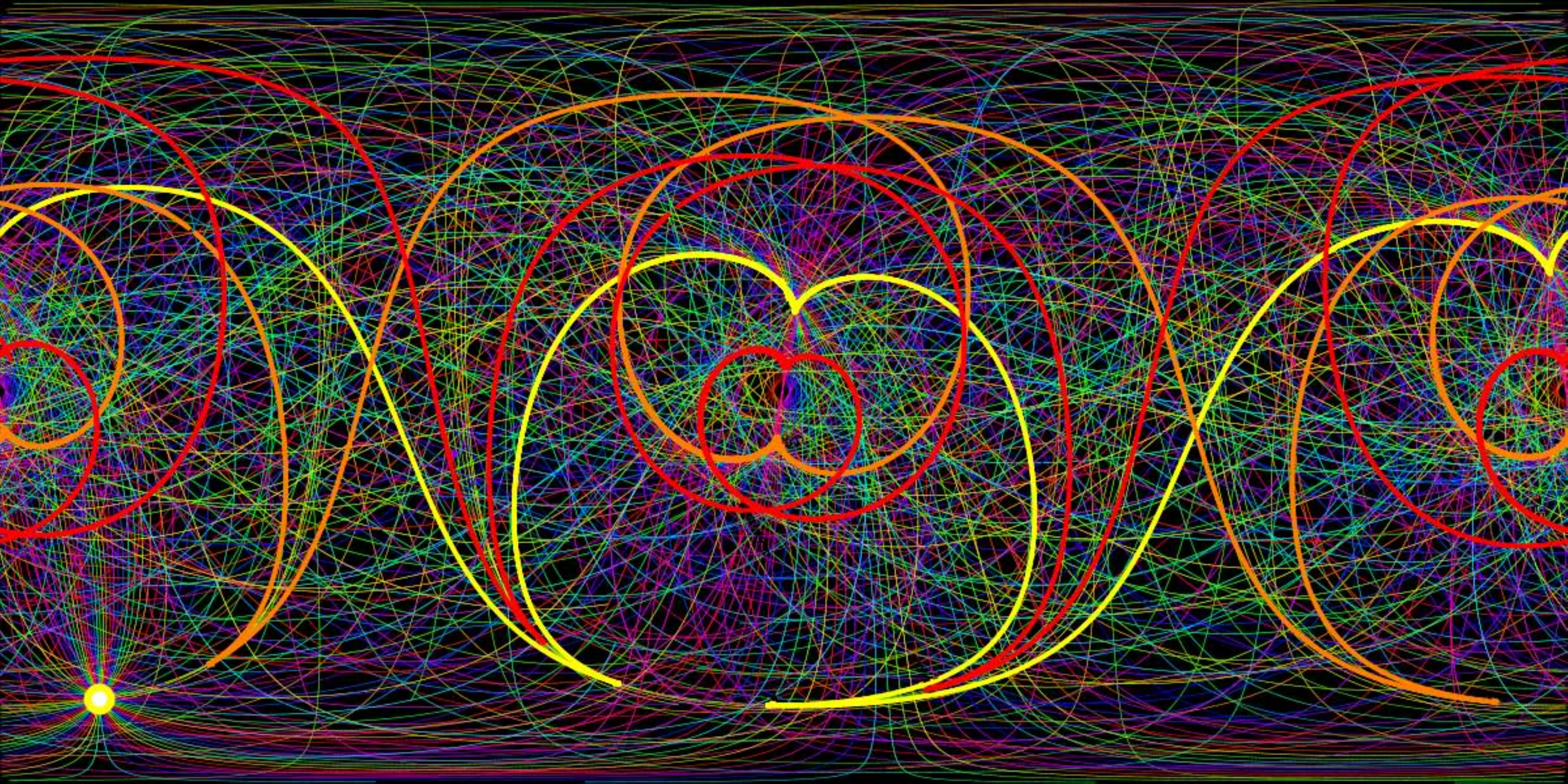}}
\caption{
An illustration to Jacobi's last geometric statement. This figure was 
computed during a HCRP project in 2009
with Michael Teodorescu. We see the first, secondary and ternary caustics of a point
on an ellipsoid $x^2/4+y^2+z^2=1$. It is an open problem whether they always have 4 cusps. In that
picture, 6000 geodesics $g(t)$ were computed solving the actual geodesic differential equation.
$g_k''= -\Gamma_k^{ij} g'_i g'_j$. The associated Gauss-Jacobi equations $f'=-K(g(t)) f$ where
$K$ is the curvature were solved too. Since special spherical coordinates were used and the
flow was forced to be on the energy surface, the numerical integration parts were done from 
scratch with Runge-Kutta rather than using built in ODE solvers. }
\end{figure}

Having a geodesic flow on surfaces or billiards allows to investigate questions in the discrete
in a completely combinatorial way which are difficult in the continuum. There are questions about
the existence of periodic orbits, the existence of geodesics which visit all points.
An other important question is how to model the continuum with discrete structures. One has
to break symmetries in order to have geometries which are not too rigid. A first possibility is
by refining randomly, an other is to use almost periodicity. The {\bf Penrose graph} for example
is a graph which is Eulerian, as all vertex degrees are 4,6,8 or 10. Taking two such Penrose patches
with smooth boundary and gluing the boundaries together produces a {\bf Penrose sphere} in $\Ecal_2$. 

\section{Platonic spheres and symmetries}

\begin{figure}[h]
\scalebox{0.16}{\includegraphics{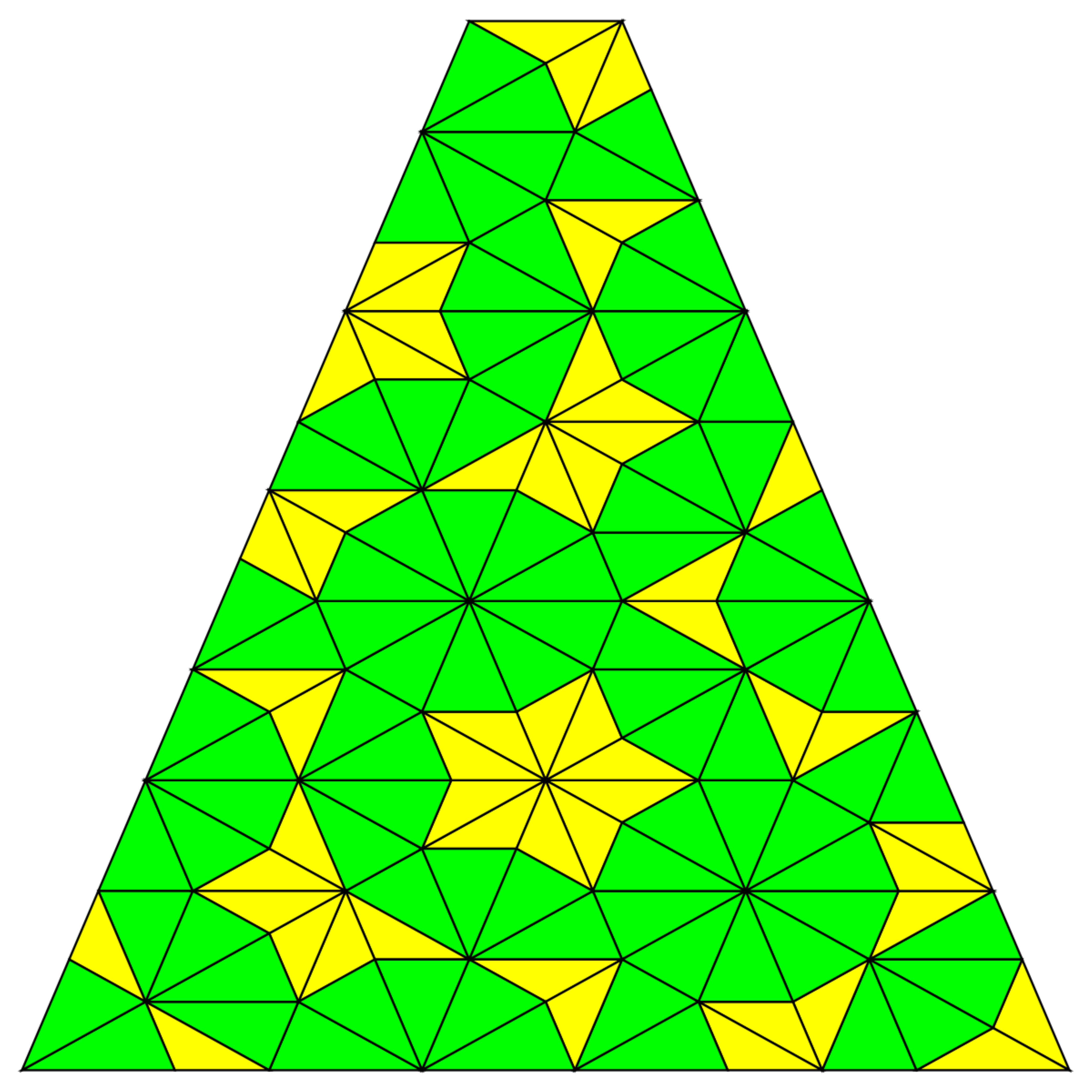}}
\caption{
A Penrose patch in which one can play billiard. The billiard is the
geodesic flow on the double cover ramified over the boundary, which is 
then an Eulerian sphere $E \in \Ecal_2$. 
}
\end{figure}

The story of regular polytopes is one of the oldest in mathematics. As pointed out
in \cite{lakatos}, the story went through an evolution of failures. The culprit is
the notion of ``polytope", a generalization of polygon and polyhedron for which is surprisingly
hard to find a precise definition. Already in school
geometry there is a mesmerizing variety of notions about what a polygon is: does it have to be 
a simple polygon or convex? For polyhedra in three dimension, one does not know whether to see 
the story as a comedy or tragedy \cite{Richeson} 
as ``theorems" like the Euler polyhedron formula had to go through a chain of improvements
when new counter examples to the formula appeared. Examples like the Kepler-Poinsot polyhedra with
negative Euler characteristic illustrate this. Of course the dust has settled today.
In higher dimensions, one usually restricts now to convex polyopes \cite{gruenbaum}. 
The story of polytopes in higher dimensions is by far settled. It is historically 
interesting that early pioneer research of polytopes in higher dimensions has been done by 
non-mainstream mathematicians like Alicia Boole Stott who had a great influence on 
Schoute \cite{Schoute} similarly than Ludwig Schl\"afli \cite{Schlafli} who influenced 
Coxeter \cite{coxeter}.  \\

When leaving the Euclidean embedding questions off the table, the story of polytopes
becomes combinatorial: polygon are $1$-spheres, polyhedra are $2$-spheres and polytopes are $d$-sphere.
Yes there are graphs in $\Gcal_2$ of higher genus which have natural immersions in $R^3$ as regular polytopes
but this is exactly what makes the classical story difficult. Lets illustrate how untangled everything
becomes when giving a graph theoretical definition of Platonic solid. 

\definition{
A graph in $\Scal_d$ is a {\bf Platonic sphere} if all 
unit spheres are isomorphic Platonic spheres in $\Scal_{d-1}$.  }

Unlike the classical classification of regular polyhedra done by mathematicians like 
Schl\"afli \cite{Schlafli} Scott, Schoute or Coxeter,
the classification of Platonic spheres in graph theory
is easy to give in an independent way. We only need Gauss-Bonnet \cite{cherngaussbonnet}:

\resultlemma{
While all graphs in $\Scal_d$ with $d \leq 1$ are Platonic, 
there are exactly two Platonic solids in $\Scal_2$ and $\Scal_3$. 
In $\Scal_d$ for $d \geq 4$ there is exactly one. 
}

\begin{proof}
The case $d=-1,0,1$ is clear. For $d=2$, the curvature has to be the 
same at every vertex and add up to $2$. It therefore better has
to be positive. As it is a fraction of $6$, it is either $1/3$
(leading to the octahedron) or $1/6$ (leading to the icosahedron). 
As  all unit spheres for a graph in $\Scal_3$ have to be regular, we know that
there are two possibilities: either the unit sphere is the octahedron, 
or the icosahedron. Both are possible, the first is the 16 cell, the second the 600 cell. 
We do not have to invoke Gauss-Bonnet as it would be useless:
the sum of curvatures would be zero as for all $G \in \Scal_3$. 
This shows inductively that in higher dimensions also,, there are 
maximally two polyhedra. But in $\Scal_4$, by Gauss-Bonnet again, the 
curvatures have to add up to $2$. Since the curvature at a vertex is
$K(x) = 1-V_0/2 + V_1/3-V_2/3+V_3/4$, where $V_k$ are the $k$-volumes
in $S(x)$. The numbers $V_k(x)$ are integers so that the curvature must
be of the form $L/12$. Indeed, for $L=1$, we can realize it with the $4$-dimensional 
cross polytope (the octahedron embedded in $R^5$), which has 
the volume data $(10, 40, 80, 80, 32)$, leading to Euler characteristic
$10-40+80-80+32=2$. If it is the unit $S(x)$ of a 5-dimensional solid,
then the curvature $K(x)$ is $1-10/2+40/3-80/4+80/5-32/6 =0$ which is
no surprise as every $5$-dimensional sphere has zero Euler characteristic. 
Can we get a 4-dimensional Platonic solid in $\Scal_4$ for which the 
unit sphere $S(x)$ is the 600-cell? Lets compute  the curvature. The 
volume data of the 600-cell are $(120,720,1200,600)$ leading to the 
curvature $1-120/2+720/3-1200/4+600/5 =1$. Since the curvature adds up
to $2$ and must be the same at every vertex, this is not possible. 
We see that we are locked in to cross polytopes from now on. 
\end{proof}

We see that when looking at Platonic solids in $\Scal_d$, we can 
explain easily why the number of Platonic solids ``thins out" in higher 
dimensions: the number of Platonic solids can not increase and Gauss-Bonnet
is the reason why we the drop from $1$ to $2$ and from $3$ to $4$. 
Since odd dimensional geometric graphs have zero Euler characteristic, these
drops could can easily take place when moving from odd to even dimensions. And
since at dimensions $4$, we are already down to one, there is nothing to do
any more. We don't miss much as classically, only the duals of the platonic solids
like the hypercube have to be placed into the picture. 
The hyper-tetrahedra $K_{d+1}$ do not count for us as Platonic solids 
because they are not spheres. Indeed, they are contractible. \\

What about semi-regular polytopes? Classically this is already a bit more
difficult to analyze as one asks that the vertex degrees are the same
everywhere and that the faces are regular polytopes. A classification
in higher dimensions has not yet been done \cite{symmetries}.  \\

Lets try to give a definition which is close to the definition of 
uniform polyhedra seen classically. Of course we can not have the same
notions as we only look at geometric graphs. Denote with ${\rm Aut}(G)$ the
{\bf automorphism group} of a graph. It consists of all graph isomorphisms
$T: G \to G$. This group plays an important role also in the continuum, when 
analyzing semi regular polytopes as it plays the role of reflection and 
rotational symmetries in the continuum. 

\definition{
A graph $G \in \Scal_d$ is called a {\bf uniform sphere} 
if ${\rm Aut}(G)$ is transitive on unit spheres in the following sense:
given two unit spheres $S(x),S(y)$, there exists $T$ such
that $T(S(x))$ intersects with $S(y)$. }

{\bf Examples.} \\
{\bf 1)} Every Platonic spheres (octahedron, icosahedron) is uniform. \\
{\bf 2)} In two dimensions, all the completions of Archimedean
solids are uniform. \\
{\bf 3)} The completions of a Catalan solid for which this completion is in $\Scal_2$
are uniform spheres.  \\

Uniform spheres can also become non-Eulerian spheres, as the completion of the
dodecahedraon or the icosahedron show.
Assume $G$ is a uniform sphere. Is the dual completion $\overline{G}$ also uniform? 
We tried to prove this for a while until we realized that the dual completion dynamics
would produce larger and larger examples of graphs with the same automorphism group
for which the transitivity of the automorphism group can not hold any more.
Classification of these graphs is not yet done. Besides prismatic families of graphs, there
should be only a finite set. \\

{\bf Examples.} \\
{\bf 1)} The dual completion $\overline{G}$ of the icosahedron $G$ is the 
"small stellated dodecahedron". In an Euclidean setting,
this solid is looked at differently as see it to be embedded in Euclidean space with 12 faces. 
We look at it as an element in $\Scal_2$: there are $v_2=12 \cdot 5=60$ faces
(each of the $12$ original faces is replaced with $5$), 
$v_0=12+20=32$ vertices (the 20 original vertices together with the 12 new centers)
and $90=30+12 \cdot 5$ edges (the 30 original edges plus 5 for each of the 
original faces). The Euler characteristic is $v_0-v_1+v_2=32-90+60=2$.
The positive curvature is located on 12 new centers where it is $1/6$. \\
{\bf 2)} The completion of the dual of the octahedron is the stellated cube, 
which is also known under the name "Tetrakishexahedron". It is an example of a
Catalan solid which happens to be in $\Ecal_2$. The curvature is zero except 
for the 6 new faces vertices of curvature $1/3$.

\begin{figure}[h]
\parbox{14cm}{
\parbox{6.2cm}{\scalebox{0.15}{\includegraphics{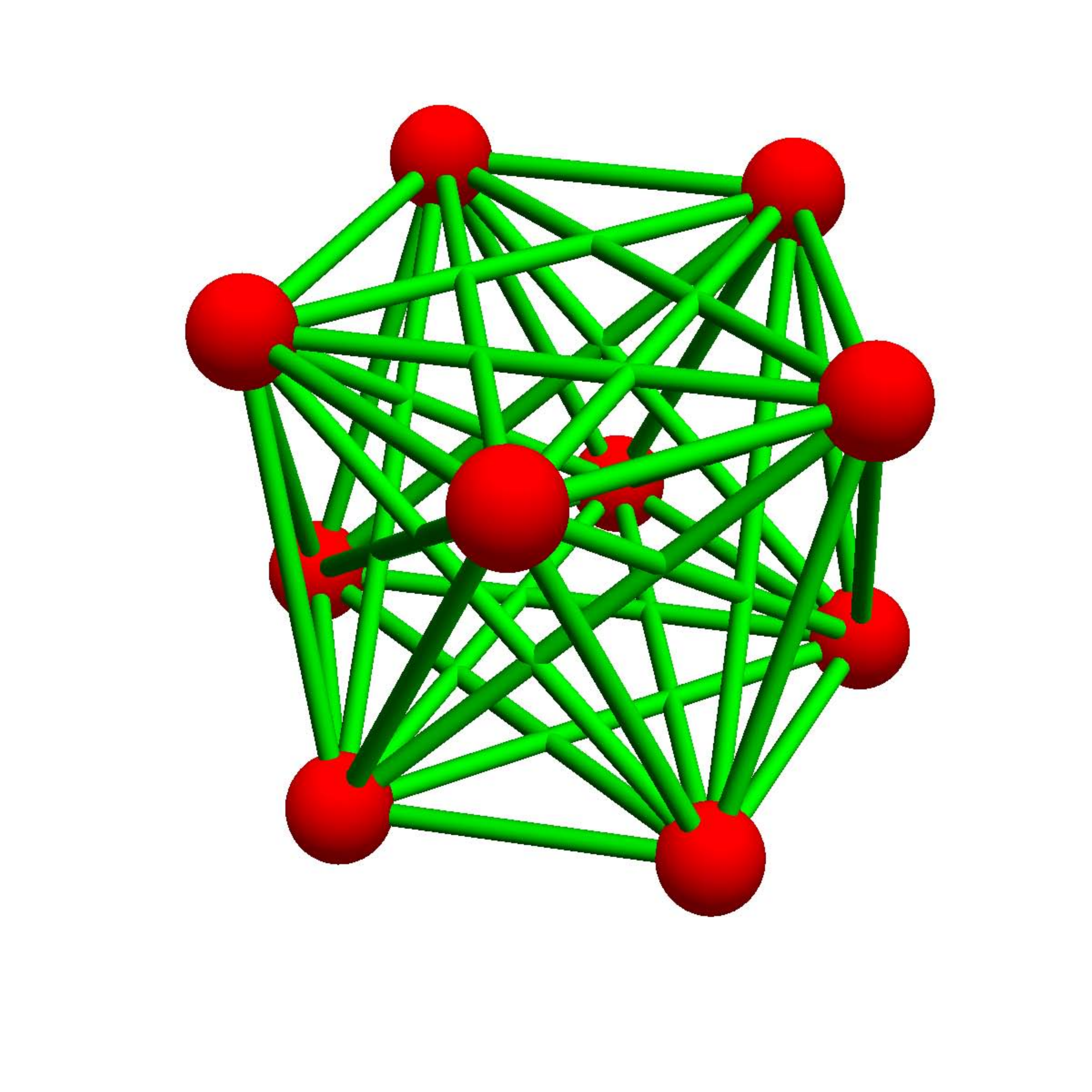}}}
\parbox{6.2cm}{\scalebox{0.17}{\includegraphics{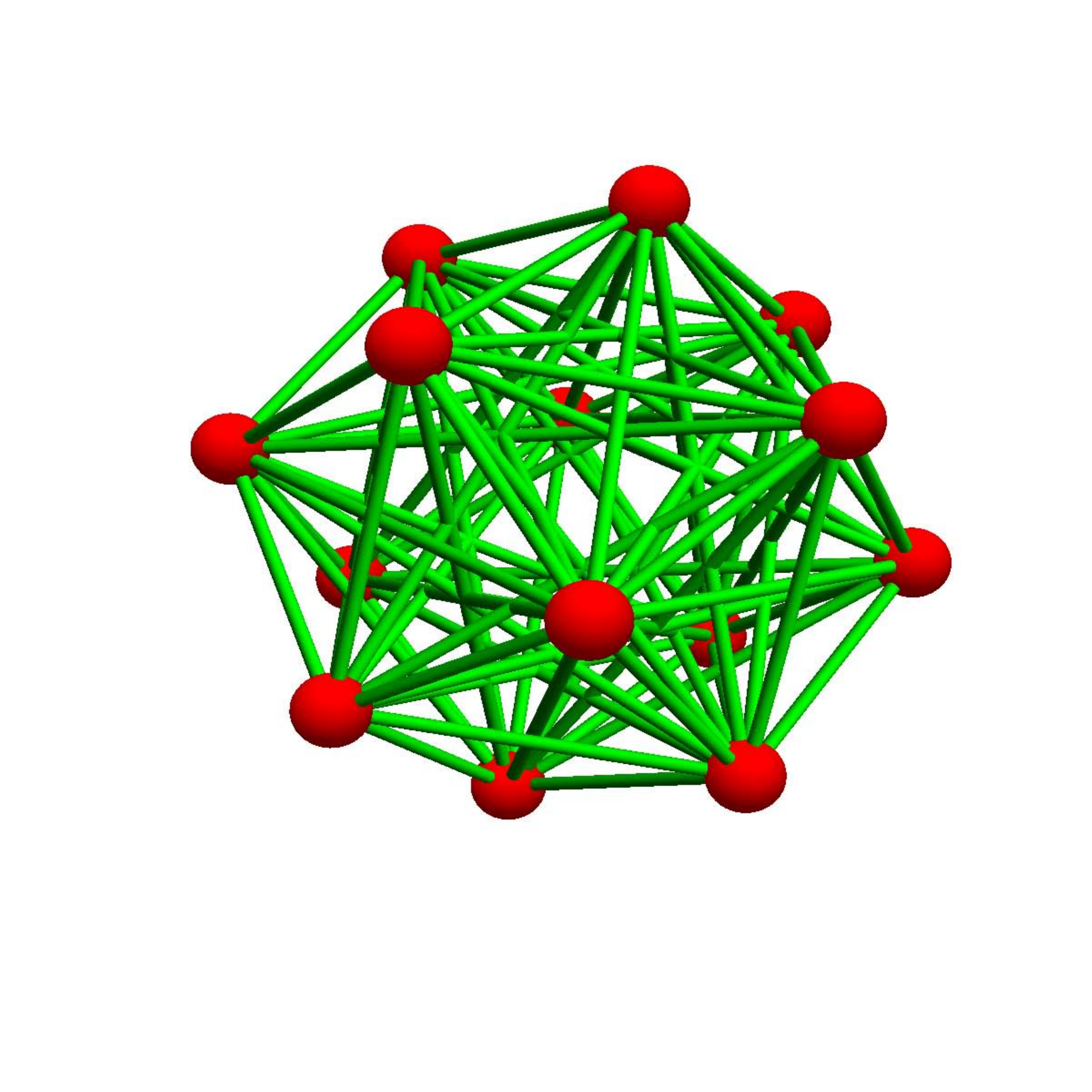}}}
}
\caption{
The 32-chamber graph in $\Scal_4$ is the only four-dimensional Platonic sphere. 
It can be obtained by drawing a regular 12-gon and connecting 
all diagonal except the 6 diagonals of maximal length. 
Once locked to a unique Platonic solid, all higher dimensional Platonic solids
are unique as the unit sphere in a Platonic solid is also Platonic, propagating
uniqueness to higher dimensions. 
To the right we see the 128-cell, the only 6-dimensional Platonic sphere. It 
can be drawn by taking a regular 14-gon and connecting all diagonals except the
main diagonals. 
}
\end{figure}

\section{Summary and conclusion}

We have looked here at the class $\Ecal_d = \Scal_d  \cap \Ccal_{d+1}$ of {\bf Eulerian spheres}, spheres which can be colored
minimally. More generally, we were interested in {\bf Eulerian geometric graphs}, geometric graphs for which all unit spheres are 
Eulerian. They are all Eulerian graphs in the classical sense that there exists a closed path visiting all edges exactly once. 
We can recursively characterize $\Ecal_d$ as the class of spheres for which all $S(x)$ are
in $\Ecal_{d-1}$ and $\Ecal_1 = \Scal_1 \cap$ $\{$ bipartite graphs $\}$ is the class of bipartite circular graphs.
We have also seen that $\Ecal_d = \Scal \cap \{ G \; | \; \hat{G} \; {\rm bipartite} \}$ 
and that $\Ecal_d \subset \Ecal \cap \Scal_d$.
We do not know yet whether $\Ecal \cap \Scal_d = \Ecal_d$. 
While we have seen that all $k$-volumes $v_k$ are even for $G \in \Scal_d$,
we do not know yet which volume data are achievable satisfying $\chi(G) = \sum_k (-1)^k v_k = 1+(-1)^d$. 
We have characterized $\Ecal_d$ as the class of spheres for which all degrees of $(d-2)$-dimensional simplices $x$ are even
where the degree is defined as the length of $\hat{x} \in \Scal_1$.
We have defined edge refinements, which are transformations from $\Scal_d \to \Scal_d$. The reverse of an edge
refinement is an edge collapse. 
An edge refinement or collapse (in the later case of course assuming that we stay in $\Gcal_d$) 
for $d \geq 2$ is a homotopy which has the effect that all degrees of maximal simplices in the $(d-2)$-dimensional 
sphere $\hat{e}$ change parity. We know by definition that any two spheres are homotopic. 
One question is whether any two spheres can be transformed into 
each other by edge refinements and collapses within $\Gcal_d$ and more generally, if any two $d$-spheres containing 
a common $(d-1)$ dimensional sphere $S$ can be transformed into each other by such transformations without
touching edges in $S$. A positive answer would verify the conjecture 
$\Scal_d \subset \Ccal_{d+2}$ which in the case $d=2$ is equivalent to the 4-color theorem.  \\

Finally, we indicated that the class of Eulerian graphs can be of interest beyond graph coloring. We characterized the
class of geometric graphs for which a Hopf-Rinov theorem holds as the class of graphs for which all unit spheres are 
projective, and especially Eulerian. In two dimensions, these graphs agree with the graphs in $\Gcal_d$ which are
Eulerian graphs in the classical sense. Having graphs with geodesic flows allows to carry virtually any question from the 
continuum to the discrete. There are questions about geodesic flows and billiards, the structure of caustics
and their relation with curvature, the size of the lengths of minimal geodesics and all within a combinatorial 
framework.  Finally we have seen that in $\Scal_d$ for $d \geq 4$, Platonic spheres, spheres for which all unit spheres
are Platonic spheres,  are unique and given by the Eulerian cross-polytopes. \\

The geometry of graphs $\Gcal_d$ and especially spheres $\Scal_d$ still needs to be investigated more. 
We have seen that the question of ``drawing lines and circles" in geometry is fundamental in order to understand
the geometry of geometric graphs, graphs for which the unit spheres are graph theoretically defined graphs. 
To ``draw lines" we need to have an exponential map which is globally defined and unique. This leads to mild
restrictions of unit spheres which are however not as severe as we just ask for a projective structure on all 
unit spheres. \\

The most prominent question remains embedding question. Can any $d$-dimensional sphere be embedded in an Eulerian 
$d+1$ dimensional sphere? Answering this question positively 
in dimension $d=2$ would prove the 4 color theorem in a geometric way. 
In higher dimensions it would lead to the conjectured bound that all $d$-spheres are either Eulerian spheres
or spheres which can be colored by $d+2$ colors.  \\

And since two dimensional manifolds can be embedded in 
the closed 4-dimensional Euclidean ball in such a way that the complement intersected with the interior 
is simply connected, one can expect that all $G \in \Gcal_2$ to have chromatic number 3,4 or 5. The 
reason is that orientable surfaces like the torus can be embedded  into the three sphere, the boundary of the 
four dimensional ball whose interior is simply connected. 
In the non-orientable case, we need the interior of a M\"obius strip to temporarily leave the boundary 3 sphere
and ``hang out" into the four dimensional interior to make a turn, but this keeps the interior simply connected. 
Coloring the four-dimensional inside with 5 colors now colors also the surface with 5 colors. 

\pagebreak

\bibliographystyle{plain}

\end{document}